\documentclass[final,3p]{elsarticle}

\usepackage{amsmath,amsthm,amssymb,mathrsfs}
\usepackage[colorlinks,linkcolor=red, citecolor=green,filecolor=magenta,urlcolor=cyan]{hyperref}

\newtheorem{theorem}{Theorem}[section]
\newtheorem{corollary}[theorem]{Corollary}
\newtheorem{remark}[theorem]{Remark}
\newtheorem{lemma}[theorem]{Lemma}
\newtheorem{proposition}[theorem]{Proposition}
\newtheorem{example}[theorem]{Example}

\numberwithin{equation}{section}

\journal{arXiv}

\allowdisplaybreaks

\begin{document}

\title{Sphere theorems for Lagrangian and Legendrian submanifolds\tnoteref{SS}}

\author[whu]{Jun Sun}
\ead{sunjun@whu.edu.cn}

\author[whu]{Linlin Sun\corref{sll}}
\ead{sunll@whu.edu.cn}

 \tnotetext[SS]{The first author was supported by the National Natural Science Foundation of China (Grant No. 11401440). Part of the work was finished when the first author was a visiting scholar at MIT supported by China Scholarship Council (CSC) and Wuhan University. The author would like to express his gratitude to Professor Tobias Colding for his invitation, to MIT for their hospitality, and to CSC and Wuhan University for their support. The second author was supported by the National Natural Science Foundation of China (Grant No. 11801420)  and Fundamental Research Funds for the Central Universities (Grant No. 2042018kf0044). }
 
 \address[whu]{School of Mathematics and Statistics \& Computational Science Hubei Key Laboratory, Wuhan University, 430072 Wuhan, China}
 
 \cortext[sll]{Corresponding author.}
 
\begin{abstract}

In this paper, we prove some differentiable sphere theorems and topological sphere theorems for Lagrangian submanifolds in K\"ahler manifold and Legendrian submanifolds in Sasaki space form.

\end{abstract}


\maketitle


\section{Introduction}

\allowdisplaybreaks

\noindent The study of Lagrangian submanifolds in a K\"ahler manifold, especially in a Calabi-Yau manifold, has attracted much attention in the past few decades (\cite{HL}, etc.), partially because of its importance in classical mechanics and mathematical physics. For instance, Strominger, Yau and Zaslow (\cite{SYZ}) found that mirror symmetry is related closely to special Lagrangian submanifolds in Calabi-Yau manifold.

\vspace{.1in}

Let $(N^{2n},\bar g,\bar \omega,J)$ be a K\"ahler manifold. A submanifold $M^n$ in $N^{2n}$ is called a {\it Lagrangian} submanifold, if the restriction of the K\"ahler form $\bar\omega$ to $M$ vanishes. Or equivalently, for any $x\in M$, $J$ maps $T_xM$ onto $N_xM$, where $J$ is the complex structure on $N$, and $T_xM$ and $N_xM$ are the tangent space and normal space of $M$ at $x$ in $N$, respectively.

Since the tangent bundle and the normal bundle of a Lagrangian submanifold are isomorphic via the complex structure J of the ambient manifold, Lagrangian submanifold has its own special properties in topology and geometry, particularly in its second fundamental form. A result of Gromov (\cite{G2}) implies that every compact embedded Lagrangian submanifold of ${\mathbb C}^n$ is not simply-connected. Of course, there exist immersed compact Lagrangian submanifolds in ${\mathbb C}^n$. One standard example is the well-known {\em Whitney sphere}, which is given by 
\begin{eqnarray*}
F:  {\mathbb S}^{n} & \longrightarrow & {\mathbb C}^n \nonumber\\
            (x_0,\cdots,x_{n})  &  \longmapsto& \frac{1}{1+x_0^2}(x_1,\cdots,x_n,x_0x_1,\cdots,x_0x_n).
\end{eqnarray*}
Gromov (\cite{G1}) also showed that a compact $n$-manifold $M$ admits a Lagrangian immersion into ${\mathbb C}^n$ if and only if the complexification of the tangent bundle, $TM\otimes {\mathbb C}$, is trivial. This can be viewed as a topological obstruction for Lagrangian submanifolds.

On the other hand, some Riemannian obstructions are also found. For example, B. Chen (\cite{Chen2}) introduced a Riemannian invariant $\delta$ in terms of the scalar curvature of $M$, and provided a sharp estimate on the invariant. A consequence of his result is that the Ricci curvature of every compact Lagrangian submanifold $M$ in ${\mathbb C}^n$ must satisfy $Ric_{\min}^M\leq 0$. Therefore, every compact irreducible symmetric space cannot be isometrically immersed in a complex Euclidean space as a Lagrangian submanifold. Furthermore, Chen showed the following sharp estimate for Lagrangian submanifolds in complex space form: 

\begin{theorem}[\cite{Chen}]
The scalar curvature $R_M$ and the mean curvature vector ${\bf H}$ of a Lagrangian submanifold in complex space form $N^{2n}$ with holomorphic sectional curvature $c$ satisfy the following sharp inequality:
\begin{equation*}
R_M\leq \frac{n(n-1)}{4}c+\frac{n-1}{n+2}|{\bf H}|^2, 
\end{equation*}
or equivalently,
\begin{equation}\label{e-lagrange}
|{\bf B}|^2\geq \frac{3}{n+2}|{\bf H}|^2,
\end{equation}
where ${\bf B}$ is the second fundamental form of $M$ in $N$
\end{theorem}

\vspace{.1in}

As we have seen, any closed Lagrangian submanifold in a complex Euclidean space cannot be isometric to a round sphere ${\mathbb S}^n$. One natural question is: when is a Lagrangian submanifold diffeomorphic, or homeomorphic to the round sphere? This involves another important topic of differential geometry: the sphere theorem. There are many interesting results on such topic (\cite{AB, BW,B, BS,BS2, CS,  GX, LS, PS, XG1,XG2}, etc.). We refer the reader to a good survey book of Brendle (\cite{B1}) for more sphere theorems under curvature pinching conditions.  In a recent paper, we also considered sphere theorems for submanifolds in K\"ahler manifold (\cite{SS}).

Based on the extra symmetries of the second fundamental form, Li-Wang (\cite{LW}) proved the following differentiable sphere theorem for Lagrangian submanifolds:

\begin{theorem}[\cite{LW}]\label{thmLW}
Let $M^n$ be an $n(\geq 3)$-dimensional compact Lagrangian submanifold in a complex space form $N^{2n}$ with holomorphic sectional curvature $c\geq 0$. Assume that
\begin{equation*}
R_M\geq\frac{(n-2)(n+1)}{4}c+\frac{2n-3}{2n+3}|{\bf H}|^2,
\end{equation*}
or equivalently,
\begin{equation*}
|{\bf B}|^2\leq \frac{3}{n+\frac{3}{2}}|{\bf H}|^2+\frac{c}{2}, 
\end{equation*}
then $M$ is diffeomorphic to a spherical space form. In particular, if $M$ is simply connected, then $M$ is diffeomorphic to ${\mathbb S}^n$.
\end{theorem}

In this paper, we wish to investigate differentiable and topological sphere theorems for Lagrangian submanifolds in K\"ahler manifold and Legendrian submanifolds in Sasaki space form under various curvature conditions.

\vspace{.1in}

Let $M$ be a smooth $n$-dimensional submanifold of a  K\"ahler manifold $N^{2m}$. We will denote the curvature tensors on $M$ and $N$ by $R$ and $K$, respectively. Recall that the sectional curvature is given by
\begin{equation*}
K(X,Y):=K(X,Y,X,Y)
\end{equation*}
and the holomorphic sectional curvature is given by
\begin{equation*}
K(X):=K(X,JX):=K(X,JX,X,JX),
\end{equation*}
where $X$ and $Y$ are tangent vector fields on $M$. Denote the minimal and maximal holomorphic sectional curvatures by
\begin{equation}\label{e1.4}
\tilde{K}_{\min}:=\min_{|X|=1}K(X),  \quad \tilde{K}_{\max}:=\max_{|X|=1}K(X).
\end{equation}
Our first theorem is about differentiable sphere theorems under Ricci curvature condition:

\vspace{.1in}

\noindent \textbf{Theorem A:} {\it Let $M$ be a smooth $n(\geq 4)$-dimensional closed simply connected Lagrangian submanifold of  a  K\"ahler manifold $N^{2n}$. If $M$ satisfies the following condition:
\begin{equation}\label{e-A}
\frac{Ric_{\min}^{[2]}}{2}\geq\dfrac{n(n-3)}{4(n-2)}\left(3\tilde K_{\max}-2\tilde K_{\min}\right)+\frac{n-3}{(n+2)(n-2)}|{\bf H}|^2,
\end{equation}
and the strict inequality holds for some point $x_0\in M$. Then $M$ is diffeomorphic to ${\mathbb S}^n$.
}

\vspace{.1in}

Here, $Ric^{[2]}$ is the 2nd weak Ricci curvature on $M$, which is a weaker assumption than the Ricci cuvature. We refer Section 2 for the definition of weak Ricci curvature.

\vspace{.1in}

By taking $\tilde{K}_{\max}=\tilde{K}_{\min}=c$, we get:
\begin{corollary}
Let $M$ be a smooth $n(\geq 4)$-dimensional closed simply connected Lagrangian submanifold of complex space form $N^{2n}$ with holomorphic sectional curvature $c$. If $M$ satisfies  the following condition:
\begin{equation*}
\frac{Ric_{\min}^{[2]}}{2}\geq\frac{n(n-3)}{4(n-2)}c+\frac{n-3}{(n+2)(n-2)}|{\bf H}|^2,
\end{equation*}
and the strict inequality holds for some point $x_0\in M$. Then $M$ is diffeomorphic to ${\mathbb S}^n$.
\end{corollary}

\vspace{.1in}

For topological sphere theorems for Lagrangian submanifolds in K\"ahler manifold, we have

\vspace{.1in}

\noindent \textbf{Theorem B:} {\it Let $M$ be a smooth $n(\geq 4)$-dimensional simply connected compact Lagrangian submanifold of  a  K\"ahler manifold $N^{2n}$. If $M$ satisfies the following condition:
\begin{equation}\label{e-B}
R_M \geq \frac{(n-3)(n+2)}{4}\left(3\tilde{K}_{\max}-2\tilde{K}_{\min}\right)+\eta(n)|{\bf H}|^2,
\end{equation}
and the strict inequality holds for some point $x_0\in M$, where $\eta(n)$ is given by
\begin{equation}\label{E-eta}
\eta(n)=
       \begin{cases}
          \hspace{0.53cm}      \frac{1}{4}, \ \ if \ n=4,\\
               \frac{3n-7}{3n+2}, \ \ if \ n\geq 5,
       \end{cases}
\end{equation} 
then $M$ is homeomorphic to ${\mathbb S}^n$.
}

\vspace{.1in}

By taking $\tilde{K}_{\max}=\tilde{K}_{\min}=c$, we obtain (comparing with Theorem \ref{thmLW}):

\begin{corollary}\label{cor1.3}
Let $M$ be a smooth $n(\geq 4)$-dimensional simply connected compact Lagrangian submanifold of complex space form $N^{2n}$ with holomorphic sectional curvature $c$. If $M$ satisfies  the following condition:
\begin{equation*}
R_M \geq \frac{(n-3)(n+2)}{4}c+\eta(n)|{\bf H}|^2,
\end{equation*}
and the strict inequality holds for some point $x_0\in M$, where $\eta(n)$ is given by (\ref{E-eta}),
then $M$ is homeomorphic to ${\mathbb S}^n$.
\end{corollary}

\begin{remark}
Recall that Hamilton (\cite{Ha}) proved that a closed simply connected four-manifold with positive isotropic curvature is diffeomorphic to ${\mathbb S}^4$. On the other hand, it is well-known that the differentiable structure on ${\mathbb S}^n$ is unique for $n=5,6$. Therefore, from the proof we see that $M$ is diffeomorphic to ${\mathbb S}^n$ for $n=4,5,6$ under the assumption of Corollary \ref{cor1.3}, which improves Theorem \ref{thmLW} for these dimensions.
\end{remark}

\vspace{.1in}

We also have topological sphere theorem under Ricci curvature conditions:

\vspace{.1in}

\noindent \textbf{Theorem C:} {\it Let $M$ be a smooth $n(\geq 4)$-dimensional simply connected compact Lagrangian submanifold of  a  K\"ahler manifold $N^{2n}$. If $M$ satisfies the following condition:
\begin{equation}\label{e-C}
Ric^{[4]}_{\min}\geq \frac{(2n-5)}{2}\left(3\tilde{K}_{\max}-2\tilde{K}_{\min}\right)+ \frac{n-3}{3n-8}|{\bf H}|^2,
\end{equation}
and the strict inequality holds for some point $x_0\in M$, 
then $M$ is homeomorphic to ${\mathbb S}^n$.
}

\begin{corollary}
Let $M$ be a smooth $n(\geq 4)$-dimensional simply connected compact Lagrangian submanifold of complex space form $N^{2n}$ with holomorphic sectional curvature $c$. If $M$ satisfies  the following condition:
\begin{equation*}
Ric^{[4]}_{\min}\geq \frac{2n-5}{2}c+\frac{n-3}{3n-8}|{\bf H}|^2,
\end{equation*}
and the strict inequality holds for some point $x_0\in M$, then $M$ is homeomorphic to ${\mathbb S}^n$.
\end{corollary}

\vspace{.1in}

\begin{remark}Assume for all orthonormal four-frames $\{e_1,e_2,e_3,e_4\}$, 
\begin{align*}
\hat{K}_{\min}\leq\dfrac14\sum_{i=1}^4K(e_i)\leq\hat{K}_{\max}.
\end{align*}
Then Theorem A, Theorem B and Theorem C are also holds if one replace $\tilde K_{\max}$ and $\tilde K_{\min}$ by $\hat{K}_{\max}$ and $\hat {K}_{\min}$ respectively.
\end{remark}

\vspace{.1in}

The main ingredient in the proof of the above theorems relies on the application of Brendle-Schoen's result on classification of closed manifold $M$ with $M\times {\mathbb R}^2$ having nonnegative isotropic curvature (\cite{BS}). In order to apply their theorem, we need to provide accurate estimates on the curvatures of the submanifolds and the ambient manifold. Similar arguments also provide sphere theorems for Legendrian submanifolds in Sasaki space forms (see Section 5).

\vspace{.1in}

The subsequent sections are organized as follows: in Section 2, we review some basic materials on K\"ahler geometry and Sasaki geometry; in Section 3, we provide the key algebraic estimates; in Section 4, we prove the sphere theorems for Lagrangian submanifolds in K\"ahler manifold; in the last Section, we obtain the sphere theorems for Legendrian submanifolds in Sasaki space form.

\vspace{.2in}

\section{Preliinaries}

\vspace{.1in}

In this section, we will provide some basic materials about K\"ahler manifold and Sasaki manifold as well as some key lemmas that will be used in the proof of the main theorems. First recall the following expression of the sectional curvature and curvature tensor in terms of holomorphic sectional curvature:

\begin{lemma}[cf. \cite{Kar}]
Let $N$ be a Riemannian manifold and $X$, $Y$, $Z$, $W$ be vector fields on $N$. Then we have
\begin{eqnarray}\label{E2.1}
24K(X,Y,Z,W)
&=& \ \ K(X+Z,Y+W)+K(X-Z,Y-W)\nonumber\\
& & +K(X+W,Y-Z)+K(X-W,Y+Z)-K(X+Z,Y-W) \nonumber\\
& & -K(X-Z,Y+W)-K(X+W,Y+Z)-K(X-W,Y-Z).
\end{eqnarray}
\end{lemma}

\begin{lemma}[cf. \cite{YK}]
Let $N$ be a K\"ahler manifold and $X$, $Y$ be vector fields on $N$. Then we have
\begin{equation}\label{E2.2}
32K(X,Y)=3K(X+JY)+3K(X-JY)-K(X+Y)-K(X-Y)-4K(X)-4K(Y).
\end{equation}
\end{lemma}

Putting (\ref{E2.2}) into (\ref{E2.1}), we get that

\begin{corollary}
Let $N$ be a K\"ahler manifold and $X, Y, Z, W$ be vector fields on $N$. Then we have
\begin{align}\label{E2.3}
256K(X,Y,Z,W)=&K(X+Z+JY+JW)+K(X+Z-JY-JW) \nonumber\\
& -K(X+Z+JY-JW)-K(X+Z-JY+JW)  \nonumber\\
& +K(X-Z+JY-JW)+K(X-Z-JY+JW) \nonumber\\
&  -K(X-Z+JY+JW)-K(X-Z-JY-JW)  \nonumber\\
& +K(X+W+JY-JZ)+K(X+W-JY+JZ) \nonumber\\
& -K(X+W+JY+JZ)-K(X+W-JY-JZ)  \nonumber\\
& +K(X-W+JY+JZ)+K(X-W-JY-JZ) \nonumber\\
&-K(X-W+JY-JZ)-K(X-W-JY+JZ).
\end{align}
\end{corollary}

\vspace{.1in}

Now we assume that $M^n$ is a Lagrangian submanifold in a K\"ahler manifold $N^{2n}$. We can choose a local orthonormal frame field of $N^{2n}$:
\begin{equation*}
\{e_1,\cdots,e_n,e_{1*},\cdots,e_{n*}\},
\end{equation*}
where $e_1,\cdots,e_n$ are tangent to $M$, $e_{1*},\cdots,e_{n*}$ are normal to $M$, and
\begin{equation*}
e_{i*}=J e_i, \ \ \ 1\leq i\leq n.
\end{equation*}
Such a frame field is called an {\it adapted frame field}. We denote
\begin{equation*}
\sigma_{ijk}:= h_{ij}^{k*}:=\langle {\bf B}(e_i,e_j),e_{k*}\rangle, 
\end{equation*}
and
\begin{equation*}
H_k:=H^{k*}:=\sum_{i=1}^n \sigma_{iik},
\end{equation*}
where $1\leq i,j,k\leq n$, and ${\bf B}$ is the second fundamental form of $M$ in $N$. Since $M$ is Lagrangian, we see that for any $X,Y,Z\in TM$
\begin{eqnarray*}
\langle {\bf B}(X,Y),JZ\rangle
&=& \langle \overline{\nabla}_XY,JZ\rangle=-\langle Y, \overline{\nabla}_X(JZ)\rangle \nonumber\\
&=& -\langle Y,J \overline{\nabla}_XZ\rangle=\langle JY,\overline{\nabla}_XZ\rangle =\langle {\bf B}(X,Z),JY\rangle.
\end{eqnarray*}
Since ${\bf B}(X,Y)$ is symmetric in $X$ and $Y$, we see that
\begin{equation*}
\langle {\bf B}(X,Y),JZ\rangle=\langle {\bf B}(Y,Z),JX\rangle=\langle {\bf B}(Z,X),JY\rangle.
\end{equation*}
In local frame, we have
\begin{equation*}
\sigma_{ijk}=\sigma_{jki}=\sigma_{kij}, \ \ \ \ 1\leq i,j,k\leq n.
\end{equation*}

\vspace{.1in}

Next we turn to Legendrian submanifolds in Sasaki space form. Let $(N^{2m+1},\phi,\xi,\eta,\bar g)$ be a compact Sasaki manifold with smooth $(1,1)$-tensor $\phi$, Reeb vector field $\xi$, contact form $\eta$ and associated Riemannian metric $\bar g$. For any unit vector field $X$ in $T_xN$ orthogonal to $\xi$, the $\phi$-sectional curvature is defined by
\begin{equation*}
H(X):=K(X,\phi X):=K(X,\phi X,X,\phi X).
\end{equation*}
A Sasaki manifold $N$ is called a {\it Sasaki space form} if $M$ has constant $\phi$-sectional curvature $c$, and will be denoted by $N(c)$. 
The following facts on Sasaki manifold will be used later. For more details, we refer the reader to the book written by Yano-Kon (\cite{YK}).

\begin{proposition}
 For a Sasaki manifold $(N^{2m+1},\phi,\xi,\eta,\bar g)$ and vector fields $X$, $Y$ on $M$, we have
\begin{align}
\eta(\xi)&=1,\notag\\
\phi^2 X&=-X+\eta(X)\xi,\notag\\
\phi\xi&=0,\notag\\
\eta(\phi X)&=0,\notag\\
\eta(X)&=\bar g(X,\xi),\notag\\
\bar g(\phi X,Y)+\bar g(X,\phi Y)&=0,\label{e2.8}\\
\overline{\nabla}_{X}\xi&=-\phi X, \label{e2.10}\\
(\overline{\nabla}_{X}\phi)Y&=\bar g(X,Y)\xi-\eta(Y)X.\label{e2.11}
\end{align}
\end{proposition}

\begin{proposition}\label{prop2.3} If the Sasaki manifold $(N^{2m+1},\phi,\xi,\eta,\bar g)$ has constant $\phi$-sectional curvature $c$, then for any vector fields $X$, $Y$, $Z$ on $M$, we have
\begin{eqnarray*}
K(X,Y)Z
&=& \frac{1}{4}(c+3)[\bar g(Y,Z)X-\bar g(X,Z)Y]\nonumber\\
& & +\frac{1}{4}(c-1)[\eta(X)\eta(Z)Y-\eta(Y)\eta(Z)X+\bar g(X,Z)\eta(Y)\xi-\bar g(Y,Z)\eta(X)\xi\nonumber\\
& & \hspace{1.7cm}+\bar g(\phi Y,Z)\phi X-\bar g(\phi X,Z)\phi Y+2\bar g(X,\phi Y)\phi Z].
\end{eqnarray*}
\end{proposition}

\vspace{.1in}

Recall that a submanifold $M^n$ in a Sasaki manifold $(N^{2n+1},\phi,\xi,\eta,\bar g)$ is called a {\it Legendrian} submanifold, if $M$ is normal to the Reeb vector field $\xi$. In particular, for any $x\in M$, $\phi$ maps $T_xM$ onto $N_xM$. For Legendrian submanifold, we have the inequality similar to (\ref{e-lagrange}) which involves the second fundamental form and mean curvature vector.

\begin{proposition}
For  an $n$-dimensional Legendrian submanifold of a $(2n+1)$-dimensional Sasaki manifold $(N^{2n+1}, \phi, \xi, \eta)$, we have
\begin{equation*}
|{\bf B}|^2\geq \frac{3}{n+2}|{\bf H}|^2.
\end{equation*}
\end{proposition}

\noindent\textbf{Proof:} Since $M^n$ is a Lengendrian submanifold in $N^{2n+1}$, we can choose an orthonormal frame field of $N^{2n+1}$:
\begin{equation*}
\{e_1,\cdots,e_n,e_{1*},\cdots,e_{n*},e_{2n+1}=\xi\},
\end{equation*}
where $e_1,\cdots,e_n$ are tangent to $M$, $e_{1*},\cdots,e_{n*},e_{2n+1}$ are normal to $M$, and
\begin{equation*}
e_{i*}=\phi e_i, \ \ \ 1\leq i\leq n.
\end{equation*}
Such a frame field is called an {\it adapted frame field}. We denote
\begin{equation*}
h_{ij}^{k*}=\langle {\bf B}(e_i,e_j),e_{k*}\rangle, \ \ \ h_{ij}^{2n+1}=\langle {\bf B}(e_i,e_j),e_{2n+1}\rangle=\langle {\bf B}(e_i,e_j),\xi\rangle,
\end{equation*}
and
\begin{equation*}
H^{k*}=\sum_{i=1}^n h_{ii}^{k*}, \ \ \ H^{2n+1}=\sum_{i=1}^n h_{ii}^{2n+1}.
\end{equation*}
Since $M$ is a Legendrian submanifold, by (\ref{e2.10}), we have for any $X,Y\in TM$ that
\begin{equation*}
\langle {\bf B}(X,Y),\xi\rangle=\langle \overline{\nabla}_XY,\xi\rangle=-\langle Y, \overline{\nabla}_X\xi\rangle=\langle Y, \phi X\rangle=0.
\end{equation*}
Hence,
\begin{equation}\label{e5.5}
h_{ij}^{2n+1}=0, \  1\leq i,j\leq n,  \ \ and \ \ H^{2n+1}=0.
\end{equation}
Similar to the Lagrangian case, set
\begin{equation*}
\sigma_{ijk}:=h_{ij}^{k*}:=\langle {\bf B}(e_i,e_j),e_{k*}\rangle, \ \ 
H_k:=H^{k*}:=\sum_{i=1}^n \sigma_{iik},
\end{equation*}
where $1\leq i,j,k\leq n$. Then (\ref{e5.5}) implies that
\begin{equation}\label{e5.9}
|{\bf B}|^2=\sum_{i,j,k=1}^n\sigma_{ijk}^2, \ \ \ |{\bf H}|^2=\sum_{k=1}^nH_k^2.
\end{equation}
Since $\eta|_M=0$, we also have, by (\ref{e2.11}) and (\ref{e2.8}), that for any $X,Y,Z\in TM$
\begin{eqnarray*}
\langle {\bf B}(X,Y),\phi Z\rangle
&=& \langle \overline{\nabla}_XY,\phi Z\rangle=-\langle Y, \overline{\nabla}_X(\phi Z)\rangle \nonumber\\
&=& -\langle Y,( \overline{\nabla}_X\phi)Z+ \phi \overline{\nabla}_XZ\rangle \nonumber\\
&=& -\langle Y,\langle X,Z\rangle\xi-\eta(Z)X\rangle-\langle Y,\phi(\nabla_XZ+{\bf B}(X,Z))\rangle\nonumber\\
&=& -\langle Y,\phi {\bf B}(X,Z)\rangle=\langle {\bf B}(X,Z),\phi Y\rangle.
\end{eqnarray*}
Since ${\bf B}(X,Y)$ is symmetric in $X$ and $Y$, we see that
\begin{equation*}
\langle {\bf B}(X,Y),\phi Z\rangle=\langle {\bf B}(Y,Z),\phi X\rangle=\langle {\bf B}(Z,X),\phi Y\rangle.
\end{equation*}
In local frame, we have
\begin{equation*}
\sigma_{ijk}=\sigma_{jki}=\sigma_{kij}, \ \ \ \ 1\leq i,j,k\leq n.
\end{equation*}
If we consider the decomposition:
\begin{equation*}
\sigma_{ijk}:=\mathring{\sigma}_{ijk}+\mu_i\delta_{jk}+\mu_j\delta_{ki}+\mu_k\delta_{ij},
\end{equation*}
where $\mu_k=\frac{1}{n+2}H_k$, then it is easy to check that 
\begin{align*}
\mathring{\sigma}_{ijk}=\mathring{\sigma}_{jik}=\mathring{\sigma}_{ikj}, \ \ \ \sum_{i=1}^n\mathring{\sigma}_{iik}=0.
\end{align*}
Moreover, (\ref{e5.9}) implies that
\begin{align*}
|{\bf B}|^2=|\sigma|^2=|\mathring{\sigma}|^2+3(n+2)|\mu|^2=|\mathring{\sigma}|^2+\frac{3}{n+2}|{\bf H}|^2,
\end{align*}
where
\begin{align*}
|\sigma|^2=\sum_{i,j,k=1}^n\sigma_{ijk}^2,\quad |\mathring{\sigma}|^2=\sum_{i,j,k=1}^n\mathring{\sigma}_{ijk}^2,\quad|\mu|^2=\sum_{k=1}^n\mu_k^2=\frac{1}{(n+2)^2}|{\bf H}|^2.
\end{align*}
This proves the proposition.
\hfill Q.E.D.

\vspace{.1in}

From the above argument, we see that in both the Lagrangian and the Legendrian cases, the Gauss equation can be written as
\begin{equation}\label{e-gauss}
R_{ijkl}=K_{ijkl}+\sum_{m=1}^n(\sigma_{ikm}\sigma_{jlm}-\sigma_{ilm}\sigma_{jkm}).
\end{equation}
In particular, the Ricci curvature and the scalar curvature satisfies
\begin{equation*}
Ric(e_i)=R_{ii}=\sum_{j=1}^nK_{ijij}+\sum_{j,k=1}^n(\sigma_{iik}\sigma_{jjk}-\sigma_{ijk}^2),
\end{equation*}
\begin{equation}\label{e-scalar}
R_M=\sum_{i,j=1}^nK_{ijij}+|{\bf H}|^2-|{\bf B}|^2.
\end{equation}

\vspace{.1in}

Fix $p\in M$, $X,Y\in T_pM$ and an orthonormal basis $\{e_1,\cdots,e_n\}$ of $T_pM$, the following notations will be used in this paper:
\begin{align*}
Ric(X,Y)=\sum_{i=1}^nR(X,e_i,Y,e_i), \ \ Ric_{jj}=Ric(e_j,e_j),\\
[e_{i_1},\cdots,e_{i_k}]=span\{e_{i_1},\cdots,e_{i_k}\}, \ \ \ \forall 1\leq i_1<i_2 < \cdots <i_k\leq n,\\
Ric^{[k]}[e_{i_1},\cdots,e_{i_k}]=\sum_{j=1}^kRic_{i_ji_j}, \ \ Ric^{[k]}_{\min}(p)=\min_{[e_{i_1},\cdots,e_{i_k}]\subset T_pM}Ric^{[k]}[e_{i_1},\cdots,e_{i_k}],
\end{align*}
where $Ric^{[k]}[e_{i_1},\cdots,e_{i_k}]$ is called the {\it $k$-th weak Ricci curvature} of $[e_{i_1},\cdots,e_{i_k}]$, which was first introduced by Gu-Xu in \cite{GX}.

\vspace{.1in}

At the end of this section, we will state some lemmas which will be crucial in the proof of our main theorems. The first result is due to Aubin:

\begin{lemma}[\cite{Au}]\label{lemma-Aubin}Let $M$ be a  compact n-dimensional Riemannian manifold. If $M$ has nonnegative Ricci curvature everywhere and has positive Ricci curvature at some point, then $M$ admits a metric with positive Ricci curvature everywhere.
\end{lemma}

A Riemannian manifold $M$ is said to have {\it nonnegative (positive, respectively) isotropic curvature}, if
\begin{equation*}
R_{1313}+R_{1414}+R_{2323}+R_{2424}-2R_{1234}\geq 0 (>0, respectively)
\end{equation*}
for all orthonormal four-frames $\{e_1,e_2.e_3.e_4\}$. This conception was introduced by Micallef-Moore and they proved that:

\begin{lemma}[\cite{MM}]\label{lemma-MM} Let $M$ be a  compact simply connected $n(\geq 4)$-dimensional Riemannian manifold which has positive isotropic curvature, then $M$ is homeomorphic to a sphere.
\end{lemma}

In addition, Micallef-Wang proved the following topological result for manifold with positive isotropic curvature:

\begin{lemma}[\cite{MW}]\label{lemma-MW} Let $M$ be a closed even-dimensional Riemannian manifold which has positive isotropic curvature, then $b_2(M)=0$.
\end{lemma}

Furthermore, Seshadri proved the following result for manifold with nonnegative isotropic curvature:

\begin{lemma}[\cite{Se}]\label{lemma-Se}Let $M$ be a compact $n$-dimensional Riemannian manifold. If $M$ has nonnegative isotropic curvature everywhere and has positive isotropic curvature at some point, then $M$ admits a metric with positive isotropic curvature.
\end{lemma}

\vspace{.1in}

The classical 1/4-differentiable sphere theorem was finally proved by Brendle-Schoen (\cite{BS}, \cite{BS2}) using Ricci flow method. They proved that

\begin{theorem}[\cite{BS}]\label{thmBS}
Let $(M,g_0)$ be a compact, locally irreducible Riemannian manifold of dimension $n(\geq 4)$ with curvature tensor $R$. Assume that $M\times{\mathbb R}^2$ has nonnegative isotropic curvature, i.e.,
\begin{equation}\label{e-BS}
R_{1313}+\lambda^2R_{1414}+\mu^2R_{2323}+\lambda^2\mu^2R_{2424}-2\lambda\mu R_{1234}\geq0
\end{equation}
for all orthonormal four-frames $\{e_1,e_2,e_3,e_4\}$ and all $\lambda, \mu\in[0,1]$. Then one of the following statements holds:

\ \ (i) $M$ is diffeomorphic to a spherical space form;

\ (ii) $n=2m$ and the universal covering of $M$ is a K\"ahler manifold biholomorphic to $\mathbb{CP}^m$;

(iii) The universal covering of $M$ is isometric to a compact symmetric space.
\end{theorem}

\vspace{.2in}

\section{Some algebraic estimates }

\vspace{.1in}

In this section, we will prove some algebraic estimates that are crucial in the proof of the main theorems.

We say that $R$ is an algebraic curvature on $\mathbb{R}^n (n\geq4)$ if $R$ is a fourth tensor such that for every $x,y,z,w\in\mathbb{R}^n$,
\begin{align*}
\begin{cases}
R(x,y,z,w)=-R(y,x,z,w)=-R(x,y,w,z)=R(z,w,x,y),\\
R(x,y,z,w)+R(y,z,x,w)+R(z,x,y,w)=0.
\end{cases}
\end{align*}

\begin{example}
 If $\sigma=(\sigma_{ijk}):\mathbb{R}^n\times\mathbb{R}^n\times\mathbb{R}^n\longrightarrow\mathbb{R}$ is a trilinear symmetric function, we obtain an algebraic curvature tensor $\tilde R$ defined by:
\begin{align*}
\tilde R_{ijkl}:=\sum_{m=1}^n\sigma_{ikm}\sigma_{jlm}-\sum_{m=1}^n\sigma_{ilm}\sigma_{jkm},\quad\forall 1\leq i, j, k, l\leq n.
\end{align*}
\end{example}

\begin{lemma}\label{lemma3.1}
Let $R$ be an algebraic curvature tensor. Suppose there is a constant $c$ such that for every four-orthonormal frame $\{e_1,e_2,e_3,e_4\}$,
\begin{align*}
R_{1212}+R_{1234}\geq c,
\end{align*}
then for every $\lambda, \mu\in[-1,1]$ and  every four-orthonormal frame $\{e_1,e_2,e_3,e_4\}$
\begin{align*}
R_{1313}+\lambda^2R_{1414}+\mu^2R_{2323}+\lambda^2\mu^2R_{2424}-2\lambda\mu R_{1234}\geq\left(1+\lambda^2\right)\left(1+\mu^2\right)c.
\end{align*}
\end{lemma}

\noindent \textbf{Proof:} A straightforward verification (we refer to \cite{SS} for a proof).
\hfill Q.E.D.

\begin{lemma}\label{lemma3.4}Let $\sigma=(\sigma_{ijk}):\mathbb{R}^n\times\mathbb{R}^n\times\mathbb{R}^n\longrightarrow\mathbb{R}$ is a trilinear symmetric function. Define $H^{i^*}:=\sum_{j=1}^n\sigma_{jji}$ and
\begin{align*}
\tilde R_{ijkl}:=\sum_{m=1}^n\sigma_{ikm}\sigma_{jlm}-\sum_{m=1}^n\sigma_{ilm}\sigma_{jkm},\quad\forall 1\leq i, j, k, l\leq n.
\end{align*}
Then for all orthonormal frame $\{e_1,e_2,e_3,e_4,\cdots, e_n\}$,
\begin{align*}
\tilde R_{1212}+\tilde R_{1234}\geq\dfrac12\left[\dfrac{6}{2n+3}\sum_{i=1}^n\left(H^{i^*}\right)^2-\sum_{i,j,k=1}^n\sigma_{ijk}^2\right].
\end{align*}
\end{lemma}

\noindent \textbf{Proof:} The proof can be found in Li-Wang's paper (\cite{LW}). For reader's convenience, we provide a new but simpler proof.

Put
\begin{align*}
\sigma_{ijk}:=\mathring{\sigma}_{ijk}+\mu_i\delta_{jk}+\mu_j\delta_{ki}+\mu_k\delta_{ij},\quad\mu_i:=\dfrac{1}{n+2}H^{i^*}.
\end{align*}
One can check that
\begin{align*}
\sum_{i,j,k=1}^n\sigma_{ijk}^2=&\sum_{i,j,k=1}^n\mathring{\sigma}_{ijk}^2+3(n+2)\sum_{i=1}^n\mu_i^2.
\end{align*}
A straightforward calculation yields
\begin{align*}
\tilde R_{1212}=&\mu_1^2+\mu_2^2+\sum_{k=1}^n\left(\mu_k^2+\left(\mathring{\sigma}_{11k}+\mathring{\sigma}_{22k}\right)\mu_k\right)+\sum_{k=1}^n\left(\mathring{\sigma}_{11k}\mathring{\sigma}_{22k}-\mathring{\sigma}_{12k}^2\right),
\end{align*}
and
\begin{align}\label{e-1234}
\tilde R_{1234}=&\sum_{k=5}^n\left(\mathring{\sigma}_{13k}\mathring{\sigma}_{24k}-\mathring{\sigma}_{14k}\mathring{\sigma}_{23k}\right)\nonumber\\
&+\left(\mathring{\sigma}_{113}-\mathring{\sigma}_{223}\right)\mathring{\sigma}_{124}-\left(\mathring{\sigma}_{114}-\mathring{\sigma}_{224}\right)\mathring{\sigma}_{123}+\left(\mathring{\sigma}_{331}-\mathring{\sigma}_{441}\right)\mathring{\sigma}_{234}-\left(\mathring{\sigma}_{332}-\mathring{\sigma}_{442}\right)\mathring{\sigma}_{134}.
\end{align}
Therefore,
\begin{align*}
\tilde R_{1212}=&\mu_1^2+\mu_2^2+\sum_{k=1}^n\left(\mu_k^2+\left(\mathring{\sigma}_{11k}+\mathring{\sigma}_{22k}\right)\mu_k\right)+\dfrac12\sum_{k=1}^n\left(\mathring{\sigma}_{11k}+\mathring{\sigma}_{22k}\right)^2\\
&-\dfrac12\sum_{i,j,k=1}^2\mathring{\sigma}_{ijk}^2-\dfrac12\sum_{i=1}^2\sum_{k=3}^n\mathring{\sigma}_{iik}^2-\sum_{k=3}^n\mathring{\sigma}_{12k}^2,\\
\tilde R_{1234}\geq&-\dfrac12\sum_{i=1}^2\sum_{j=3}^4\sum_{k=5}^n\mathring{\sigma}_{ijk}^2-2\left(\sum_{j=3}^4\mathring{\sigma}_{12j}^2+\sum_{i=1}^2\mathring{\sigma}_{i34}^2\right)-\dfrac18\left[\sum_{j=3}^4\left(\mathring{\sigma}_{11j}-\mathring{\sigma}_{22j}\right)^2+\sum_{i=1}^2\left(\mathring{\sigma}_{i33}-\mathring{\sigma}_{i44}\right)^2\right].
\end{align*}
We obtain
\begin{align*}
\tilde R_{1212}+\tilde R_{1234}\geq&\mu_1^2+\mu_2^2+\sum_{k=1}^n\left(\mu_k^2+\left(\mathring{\sigma}_{11k}+\mathring{\sigma}_{22k}\right)\mu_k\right)+\dfrac12\sum_{k=1}^n\left(\mathring{\sigma}_{11k}+\mathring{\sigma}_{22k}\right)^2\\
&-\dfrac12\sum_{i,j,k=1}^2\mathring{\sigma}_{ijk}^2-\dfrac12\sum_{i=1}^2\sum_{k=3}^n\mathring{\sigma}_{iik}^2-3\sum_{1\leq i<j<k\leq n}\mathring{\sigma}_{ijk}^2\\
&-\dfrac18\left[\sum_{j=3}^4\left(\mathring{\sigma}_{11j}-\mathring{\sigma}_{22j}\right)^2+\sum_{i=1}^2\left(\mathring{\sigma}_{i33}-\mathring{\sigma}_{i44}\right)^2\right]\\
\geq&\mu_1^2+\mu_2^2+\sum_{k=1}^n\left(\mu_k^2+\left(\mathring{\sigma}_{11k}+\mathring{\sigma}_{22k}\right)\mu_k\right)+\dfrac12\sum_{k=1}^n\left(\mathring{\sigma}_{11k}+\mathring{\sigma}_{22k}\right)^2\\
&-\dfrac12\sum_{i,j,k=1}^n\mathring{\sigma}_{ijk}^2+\dfrac12\left(\sum_{i=3}^n\mathring{\sigma}_{iii}^2+3\sum_{3\leq i\neq j\leq n}\mathring{\sigma}_{ijj}^2\right)+\sum_{i=1}^2\sum_{k=3}^n\mathring{\sigma}_{iik}^2\\
&+\dfrac32\sum_{i=1}^2\sum_{k=3}^n\mathring{\sigma}_{ikk}^2-\dfrac18\left[\sum_{j=3}^4\left(\mathring{\sigma}_{11j}-\mathring{\sigma}_{22j}\right)^2+\sum_{i=1}^2\left(\mathring{\sigma}_{i33}-\mathring{\sigma}_{i44}\right)^2\right]\\
\geq&\mu_1^2+\mu_2^2+\sum_{k=1}^n\left(\mu_k^2+\left(\mathring{\sigma}_{11k}+\mathring{\sigma}_{22k}\right)\mu_k\right)+\dfrac12\sum_{k=1}^n\left(\mathring{\sigma}_{11k}+\mathring{\sigma}_{22k}\right)^2\\
&-\dfrac12\sum_{i,j,k=1}^n\mathring{\sigma}_{ijk}^2+\dfrac{3}{2n}\sum_{i=3}^n\left(\sum_{j=3}^n\mathring{\sigma}_{ijj}\right)^2+\dfrac12\sum_{k=3}^n\left(\sum_{i=1}^2\mathring{\sigma}_{iik}\right)^2\\
&+\dfrac34\sum_{i=1}^2\left(\sum_{k=3}^4\mathring{\sigma}_{ikk}\right)^2+\dfrac32\sum_{i=1}^2\sum_{k=5}^n\mathring{\sigma}_{ikk}^2\\
\geq&\mu_1^2+\mu_2^2+\sum_{k=1}^n\left(\mu_k^2+\left(\mathring{\sigma}_{11k}+\mathring{\sigma}_{22k}\right)\mu_k\right)\\
&-\dfrac12\sum_{i,j,k=1}^n\mathring{\sigma}_{ijk}^2+\dfrac{2n+3}{2n}\sum_{i=3}^n\left(\sum_{j=1}^2\mathring{\sigma}_{ijj}\right)^2+\dfrac{n+1}{2(n-2)}\sum_{i=1}^2\left(\sum_{k=3}^n\mathring{\sigma}_{ikk}\right)^2\\
\geq&\dfrac{3(n+2)}{2(n+1)}\sum_{k=1}^2\mu_k^2+\dfrac{3(n+2)}{2(2n+3)}\sum_{k=3}^n\mu_k^2-\dfrac12\sum_{i,j,k=1}^n\mathring{\sigma}_{ijk}^2.
\end{align*}
We conclude that
\begin{align*}
\tilde R_{1212}+\tilde R_{1234}\geq&\dfrac{3(n+2)}{2(2n+3)}\sum_{k=1}^n\mu_k^2-\dfrac12\sum_{i,j,k=1}^n\mathring{\sigma}_{ijk}^2\\
=&\dfrac{3}{2n+3}\sum_{k=1}^n\left(H^{k^*}\right)^2-\dfrac12\sum_{i,j,k=1}^n\sigma_{ijk}^2.
\end{align*}
\hfill Q.E.D.

\begin{lemma}\label{lemma3.5}Let $\sigma$ and $\tilde R$ be as in Lemma \ref{lemma3.4}. Then for all orthonormal frame $\{e_1,e_2,e_3,e_4,\cdots, e_n\}$,
\begin{align*}
\sum_{i=1}^2\sum_{j=3}^4\tilde R_{ijij}-2\tilde R_{1234}\geq&\tilde\eta(n)\sum_{k=1}^n\left(H^{k^*}\right)^2-\dfrac{2}{3}\sum_{i,j,k=1}^n\sigma_{ijk}^2,
\end{align*}
where
\begin{equation*}
\tilde\eta(n)=
       \begin{cases}
          \hspace{0.53cm}      \frac{1}{2}, \ \ if \ n=4,\\
               \frac{6}{3n+2}, \ \ if \ n\geq 5.
       \end{cases}
\end{equation*} 
\end{lemma}

\noindent \textbf{Proof:} Using the same notations as in the proof of Lemma \ref{lemma3.4}, a direct calculation yields
\begin{align*}
\sum_{i,j,k=1}^n\mathring{\sigma}_{ijk}^2=
&\left(\sum_{i=1}^2\mathring{\sigma}_{iii}^2+3\sum_{1\leq i\neq j\leq 2}\mathring{\sigma}_{iij}^2\right)+\left(\sum_{i=1}^3\mathring{\sigma}_{iii}^4+3\sum_{3\leq i\neq j\leq 4}\mathring{\sigma}_{iij}^2\right)\\
&+\left(\sum_{i=5}^n\mathring{\sigma}_{iii}^2+3\sum_{5\leq i\neq j\leq n}\mathring{\sigma}_{iij}^2\right)+3\sum_{i=1}^2\sum_{j=3}^4\mathring{\sigma}_{iij}^2+3\sum_{i=1}^2\sum_{j=3}^4\mathring{\sigma}_{ijj}^2\\
&+3\sum_{i=1}^4\sum_{j=5}^n\mathring{\sigma}_{iij}^2+3\sum_{i=1}^4\sum_{j=5}^n\mathring{\sigma}_{ijj}^2+6\sum_{1\leq i<j<k\leq n}\mathring{\sigma}_{ijk}^2.
\end{align*}
On one hand, notice that
\begin{align}\label{e-cur}
\sum_{i=1}^2\sum_{j=3}^4\tilde{R}_{ijij}=&2\sum_{i=1}^4\mu_i^2+4\sum_{k=1}^n\mu_k^2+2\sum_{k=1}^n\sum_{i=1}^4\mathring{\sigma}_{iik}\mu_k+\sum_{i=1}^2\sum_{j=3}^4\sum_{k=1}^n\left(\mathring{\sigma}_{iik}\mathring{\sigma}_{jjk}-\mathring{\sigma}_{ijk}^2\right)\nonumber\\
=&2\sum_{i=1}^4\mu_i^2+4\sum_{k=1}^n\mu_k^2+2\sum_{k=1}^n\sum_{i=1}^4\mathring{\sigma}_{iik}\mu_k+\dfrac12\sum_{k=1}^n\left(\sum_{i=1}^4\mathring{\sigma}_{iik}\right)^2\nonumber\\
&-\dfrac12\sum_{k=1}^n\left(\mathring{\sigma}_{11k}+\mathring{\sigma}_{22k}\right)^2-\dfrac12\sum_{k=1}^n\left(\mathring{\sigma}_{33k}+\mathring{\sigma}_{44k}\right)^2-\sum_{i=1}^2\sum_{j=3}^4\sum_{k=1}^n\mathring{\sigma}_{ijk}^2\nonumber\\
=&2\sum_{i=1}^4\mu_i^2+4\sum_{k=1}^n\mu_k^2+2\sum_{k=1}^n\sum_{i=1}^4\mathring{\sigma}_{iik}\mu_k+\dfrac12\sum_{k=1}^n\left(\sum_{i=1}^4\mathring{\sigma}_{iik}\right)^2\nonumber\\
&-\dfrac12\sum_{k=1}^2\left(\mathring{\sigma}_{11k}+\mathring{\sigma}_{22k}\right)^2-\dfrac12\sum_{k=3}^4\left(\mathring{\sigma}_{11k}+\mathring{\sigma}_{22k}\right)^2-\dfrac12\sum_{k=5}^n\left(\mathring{\sigma}_{11k}+\mathring{\sigma}_{22k}\right)^2\nonumber\\
&-\dfrac12\sum_{k=1}^2\left(\mathring{\sigma}_{33k}+\mathring{\sigma}_{44k}\right)^2-\dfrac12\sum_{k=3}^4\left(\mathring{\sigma}_{33k}+\mathring{\sigma}_{44k}\right)^2-\dfrac12\sum_{k=5}^n\left(\mathring{\sigma}_{33k}+\mathring{\sigma}_{44k}\right)^2\nonumber\\
&-\sum_{i=1}^2\sum_{j=3}^4\mathring{\sigma}_{iij}^2-\sum_{i=1}^2\sum_{j=3}^4\mathring{\sigma}_{ijj}^2-\sum_{i=1}^2\sum_{j=3}^4\sum_{k=5}^n\mathring{\sigma}_{ijk}^2-2\sum_{k=3}^4\mathring{\sigma}_{12k}^2-2\sum_{i=1}^2\mathring{\sigma}_{i34}^2\nonumber\\
=&2\sum_{i=1}^4\mu_i^2+4\sum_{k=1}^n\mu_k^2+2\sum_{k=1}^n\sum_{i=1}^4\mathring{\sigma}_{iik}\mu_k+\dfrac12\sum_{k=1}^n\left(\sum_{i=1}^4\mathring{\sigma}_{iik}\right)^2\nonumber\\
&-\sum_{k=1}^2\left(\mathring{\sigma}_{33k}+\mathring{\sigma}_{44k}\right)^2-\sum_{k=3}^4\left(\mathring{\sigma}_{11k}+\mathring{\sigma}_{22k}\right)^2-\dfrac12\sum_{k=5}^n\left(\mathring{\sigma}_{11k}+\mathring{\sigma}_{22k}\right)^2\nonumber\\
&-\dfrac12\sum_{k=1}^2\left(\mathring{\sigma}_{11k}+\mathring{\sigma}_{22k}\right)^2-\dfrac12\sum_{k=3}^4\left(\mathring{\sigma}_{33k}+\mathring{\sigma}_{44k}\right)^2-\dfrac12\sum_{k=5}^n\left(\mathring{\sigma}_{33k}+\mathring{\sigma}_{44k}\right)^2\nonumber\\
&-\dfrac12\sum_{k=3}^4\left(\mathring{\sigma}_{11k}-\mathring{\sigma}_{22k}\right)^2-\dfrac12\sum_{k=1}^2\left(\mathring{\sigma}_{33k}-\mathring{\sigma}_{44k}\right)^2\nonumber\\
& -2\sum_{k=3}^4\mathring{\sigma}_{12k}^2-2\sum_{i=1}^2\mathring{\sigma}_{i34}^2-\sum_{i=1}^2\sum_{j=3}^4\sum_{k=5}^n\mathring{\sigma}_{ijk}^2.
\end{align}
But,
\begin{align*}
\tilde{R}_{1234}=&\sum_{k=5}^n\left(\mathring{\sigma}_{13k}\mathring{\sigma}_{24k}-\mathring{\sigma}_{14k}\mathring{\sigma}_{23k}\right)\\
&+\left(\mathring{\sigma}_{113}-\mathring{\sigma}_{223}\right)\mathring{\sigma}_{124}-\left(\mathring{\sigma}_{114}-\mathring{\sigma}_{224}\right)\mathring{\sigma}_{123}+\left(\mathring{\sigma}_{331}-\mathring{\sigma}_{441}\right)\mathring{\sigma}_{234}-\left(\mathring{\sigma}_{332}-\mathring{\sigma}_{442}\right)\mathring{\sigma}_{134}\\
\leq&\left(\sum_{j=3}^4\mathring{\sigma}_{12j}^2+\sum_{i=1}^2\mathring{\sigma}_{i34}^2\right)+\dfrac{1}{4}\left(\sum_{j=3}^4\left(\mathring{\sigma}_{11j}-\mathring{\sigma}_{22j}\right)^2+\sum_{i=1}^2\left(\mathring{\sigma}_{33i}-\mathring{\sigma}_{44i}\right)^2\right)\\
& +\dfrac12\sum_{i=1}^2\sum_{j=3}^4\sum_{k=5}^n\mathring{\sigma}_{ijk}^2.
\end{align*}
We obtain
\begin{align*}
\sum_{i=1}^2\sum_{j=3}^4\tilde{R}_{ijij}-2\tilde{R}_{1234}\geq&2\sum_{i=1}^4\mu_i^2+4\sum_{k=1}^n\mu_k^2+2\sum_{k=1}^n\sum_{i=1}^4\mathring{\sigma}_{iik}\mu_k+\dfrac12\sum_{k=1}^n\left(\sum_{i=1}^4\mathring{\sigma}_{iik}\right)^2\\
&-\sum_{k=1}^2\left(\mathring{\sigma}_{33k}+\mathring{\sigma}_{44k}\right)^2-\sum_{k=3}^4\left(\mathring{\sigma}_{11k}+\mathring{\sigma}_{22k}\right)^2-\dfrac12\sum_{k=5}^n\left(\mathring{\sigma}_{11k}+\mathring{\sigma}_{22k}\right)^2\\
&-\dfrac12\sum_{k=1}^2\left(\mathring{\sigma}_{11k}+\mathring{\sigma}_{22k}\right)^2-\dfrac12\sum_{k=3}^4\left(\mathring{\sigma}_{33k}+\mathring{\sigma}_{44k}\right)^2-\dfrac12\sum_{k=5}^n\left(\mathring{\sigma}_{33k}+\mathring{\sigma}_{44k}\right)^2\\
&-\sum_{k=3}^4\left(\mathring{\sigma}_{11k}-\mathring{\sigma}_{22k}\right)^2-\sum_{k=1}^2\left(\mathring{\sigma}_{33k}-\mathring{\sigma}_{44k}\right)^2-4\sum_{1\leq i<j<k\leq n}\mathring{\sigma}_{ijk}^2.
\end{align*}
On the other hand,
\begin{align*}
& \dfrac23\sum_{i,j,k=1}^n\mathring{\sigma}_{ijk}^2\\
=&\dfrac23\left(\sum_{i=1}^2\mathring{\sigma}_{iii}^2+3\sum_{1\leq i\neq j\leq 2}\mathring{\sigma}_{iij}^2\right)+\dfrac23\left(\sum_{i=3}^4\mathring{\sigma}_{iii}^2+3\sum_{3\leq i\neq j\leq 4}\mathring{\sigma}_{iij}^2\right)+\dfrac23\left(\sum_{i=5}^n\mathring{\sigma}_{iii}^2+3\sum_{5\leq i\neq j\leq n}\mathring{\sigma}_{iij}^2\right)\\
&+\sum_{k=3}^4\left(\mathring{\sigma}_{11k}+\mathring{\sigma}_{22k}\right)^2+\sum_{k=1}^2\left(\mathring{\sigma}_{33k}+\mathring{\sigma}_{44k}\right)^2+2\sum_{i=1}^4\sum_{j=5}^n\mathring{\sigma}_{iij}^2+2\sum_{i=1}^4\sum_{j=5}^n\mathring{\sigma}_{ijj}^2\\
&+\sum_{k=3}^4\left(\mathring{\sigma}_{11k}-\mathring{\sigma}_{22k}\right)^2+\sum_{k=1}^2\left(\mathring{\sigma}_{33k}-\mathring{\sigma}_{44k}\right)^2+4\sum_{1\leq i<j<k\leq n}\mathring{\sigma}_{ijk}^2\\
\geq&\dfrac12\sum_{i=1}^2\left(\mathring{\sigma}_{11i}+\mathring{\sigma}_{22i}\right)^2+\dfrac12\sum_{i=3}^4\left(\mathring{\sigma}_{33i}+\mathring{\sigma}_{44i}\right)^2+\dfrac23\left(\sum_{i=5}^n\mathring{\sigma}_{iii}^2+3\sum_{5\leq i\neq j\leq n}\mathring{\sigma}_{iij}^2\right)\\
&+\sum_{k=3}^4\left(\mathring{\sigma}_{11k}+\mathring{\sigma}_{22k}\right)^2+\sum_{k=1}^2\left(\mathring{\sigma}_{33k}+\mathring{\sigma}_{44k}\right)^2+\sum_{i=1}^4\sum_{j=5}^n\mathring{\sigma}_{iij}^2+2\sum_{i=1}^4\sum_{j=5}^n\mathring{\sigma}_{ijj}^2\\
&+\dfrac12\sum_{j=5}^n\left[\left(\mathring{\sigma}_{11j}+\mathring{\sigma}_{22j}\right)^2+\left(\mathring{\sigma}_{33j}+\mathring{\sigma}_{44j}\right)^2\right]+\sum_{k=3}^4\left(\mathring{\sigma}_{11k}-\mathring{\sigma}_{22k}\right)^2+\sum_{k=1}^2\left(\mathring{\sigma}_{33k}-\mathring{\sigma}_{44k}\right)^2\\
& +4\sum_{1\leq i<j<k\leq n}\mathring{\sigma}_{ijk}^2.
\end{align*}

\vspace{.1in}

If $n=4$, we have 
\begin{equation*}
\sum_{i=1}^4\mathring{\sigma}_{iik}=0, \ \ \ for \ 1\leq k\leq 4.
\end{equation*}
Then we compute
\begin{align*}
\sum_{i=1}^2\sum_{j=3}^4\tilde{R}_{ijij}-2\tilde{R}_{1234}
\geq&6\sum_{k=1}^4\mu_k^2
-\sum_{k=1}^2\left(\mathring{\sigma}_{33k}+\mathring{\sigma}_{44k}\right)^2-\sum_{k=3}^4\left(\mathring{\sigma}_{11k}+\mathring{\sigma}_{22k}\right)^2\\
&-\dfrac12\sum_{k=1}^2\left(\mathring{\sigma}_{11k}+\mathring{\sigma}_{22k}\right)^2-\dfrac12\sum_{k=3}^4\left(\mathring{\sigma}_{33k}+\mathring{\sigma}_{44k}\right)^2\\
&-\sum_{k=3}^4\left(\mathring{\sigma}_{11k}-\mathring{\sigma}_{22k}\right)^2-\sum_{k=1}^2\left(\mathring{\sigma}_{33k}-\mathring{\sigma}_{44k}\right)^2-4\sum_{1\leq i<j<k\leq 4}\mathring{\sigma}_{ijk}^2\\
\geq & 6\sum_{k=1}^4\mu_k^2-\frac{2}{3}\sum_{i,j,k=1}^4\mathring{\sigma}_{ijk}^2\\
= &\dfrac{1}{2}\sum_{k=1}^4\left(H^{k^*}\right)^2-\dfrac{2}{3}\sum_{i,j,k=1}^n\sigma_{ijk}^2.
\end{align*}

\vspace{.1in}

If $n\geq 5$, we compute
\begin{align*}
&\sum_{i=1}^2\sum_{j=3}^4\tilde{R}_{ijij}-2\tilde{R}_{1234}+\dfrac23\sum_{i,j,k=1}^n\mathring{\sigma}_{ijk}^2\\
\geq&\dfrac12\sum_{k=1}^n\left(\sum_{i=1}^4\mathring{\sigma}_{iik}\right)^2+2\sum_{k=1}^n\sum_{i=1}^4\mathring{\sigma}_{iik}\mu_k+2\sum_{i=1}^4\mu_i^2+4\sum_{k=1}^n\mu_k^2\\
&+\dfrac23\left(\sum_{i=5}^n\mathring{\sigma}_{iii}^2+3\sum_{5\leq i\neq j\leq n}\mathring{\sigma}_{iij}^2\right)+\sum_{i=1}^4\sum_{j=5}^n\mathring{\sigma}_{iij}^2+2\sum_{i=1}^4\sum_{j=5}^n\mathring{\sigma}_{ijj}^2\\
\geq&\dfrac12\sum_{k=1}^n\left(\sum_{i=1}^4\mathring{\sigma}_{iik}\right)^2+2\sum_{k=1}^n\sum_{i=1}^4\mathring{\sigma}_{iik}\mu_k+2\sum_{i=1}^4\mu_i^2+4\sum_{k=1}^n\mu_k^2\\
&+\dfrac{2}{n-2}\sum_{j=5}^n\left(\sum_{k=5}^n\mathring{\sigma}_{jkk}\right)^2+\dfrac14\sum_{j=5}^n\left(\sum_{i=1}^4\mathring{\sigma}_{iij}\right)^2+\dfrac{2}{n-4}\sum_{i=1}^4\left(\sum_{j=5}^n\mathring{\sigma}_{ijj}\right)^2\\
=&2\sum_{k=1}^n\sum_{i=1}^4\mathring{\sigma}_{iik}\mu_k+2\sum_{i=1}^4\mu_i^2+4\sum_{k=1}^n\mu_k^2\\
&+\dfrac{3n+2}{4(n-2)}\sum_{k=5}^n\left(\sum_{i=1}^4\mathring{\sigma}_{iik}\right)^2+\dfrac{n}{2(n-4)}\sum_{k=1}^4\left(\sum_{i=1}^4\mathring{\sigma}_{iik}\right)^2\\
\geq&4\sum_{k=1}^n\mu_k^2+\dfrac{8}{n}\sum_{i=1}^4\mu_i^2-\dfrac{4(n-2)}{3n+2}\sum_{j=5}^n\mu_j^2\\
\geq&\dfrac{8(n+2)}{3n+2}\sum_{k=1}^n\mu_k^2.
\end{align*}
Finally, we obtain
\begin{align*}
\sum_{i=1}^2\sum_{j=3}^4\tilde{R}_{ijij}-2\tilde{R}_{1234}\geq\dfrac{6}{3n+2}\sum_{k=1}^n\left(H^{k^*}\right)^2-\dfrac{2}{3}\sum_{i,j,k=1}^n\sigma_{ijk}^2.
\end{align*}
This finishes the proof of the lemma.
\hfill Q.E.D.

\begin{lemma}\label{lemma3.6}Let $\sigma$ and $\tilde R$ be as in Lemma \ref{lemma3.4}. Then for for all orthonormal frame $\{e_1,e_2,e_3,e_4,\cdots, e_n\}$,
\begin{align*}
\tilde R_{1212}+\tilde R_{1234}\geq\dfrac{1}{2}\sum_{i=1}^2\tilde R_{ii}+\dfrac{n-4}{2(n-2)}\sum_{i=3}^n\tilde R_{ii}-\dfrac{n-3}{2(n+2)}\sum_{k=1}^n\left(H^{k^*}\right)^2,
\end{align*}
where $\tilde R_{ii}:=\sum_{k=1}^n\tilde R_{ikik}$.
\end{lemma}

\noindent \textbf{Proof:} As notations in Lemma \ref{lemma3.4}, we get
\begin{align*}
\tilde R_{ii}:=\sum_{k=1}^n\tilde R_{ikik}=(n-2)\mu_i^2+(n-2)\sum_{k=1}^n\mathring{\sigma}_{iik}\mu_k+n\sum_{k=1}^n\mu_k^2-\sum_{j,k=1}^n\mathring{\sigma}_{ijk}^2.
\end{align*}
Thus,
\begin{align*}
\dfrac{1}{2}\sum_{i=1}^2\tilde R_{ii}=&\dfrac{n-2}{2}\left(\sum_{k=1}^n\sum_{i=1}^2\mathring{\sigma}_{iik}\mu_k+\sum_{i=1}^2\mu_i^2\right)+n\sum_{k=1}^n\mu_k^2-\dfrac12\sum_{i=1}^2\sum_{j,k=1}^n\mathring{\sigma}_{ijk}^2,\\
\dfrac{1}{n-2}\sum_{i=3}^n\tilde R_{ii}=&\left(\sum_{k=1}^n\sum_{i=3}^n\mathring{\sigma}_{iik}\mu_k+\sum_{i=3}^n\mu_i^2\right)+n\sum_{k=1}^n\mu_k^2-\dfrac{1}{n-2}\sum_{i=3}^n\sum_{j,k=1}^n\mathring{\sigma}_{ijk}^2.
\end{align*}
By assumption,
\begin{align*}
&\dfrac{\varepsilon}{2}\sum_{i=1}^2\tilde R_{ii}+\dfrac{1-\varepsilon}{n-2}\sum_{i=3}^n\tilde R_{ii}\\
=&\dfrac{n\varepsilon-2}{2}\sum_{k=1}^n\sum_{i=1}^2\mathring{\sigma}_{iik}\mu_k+\dfrac{(n-2)\varepsilon}{2}\sum_{i=1}^2\mu_i^2+(1-\varepsilon)\sum_{i=3}^n\mu_i^2+n\sum_{k=1}^n\mu_k^2\\
&-\dfrac{\varepsilon}{2}\sum_{i=1}^2\sum_{j,k=1}^n\mathring{\sigma}_{ijk}^2-\dfrac{1-\varepsilon}{n-2}\sum_{i=3}^n\sum_{j,k=1}^n\mathring{\sigma}_{ijk}^2\\
=&\dfrac{(n-2)\varepsilon-2}{2}\sum_{k=1}^n\sum_{i=1}^2\mathring{\sigma}_{iik}\mu_k+\dfrac{(n-4)\varepsilon}{2}\sum_{i=1}^2\mu_i^2+(1-\varepsilon)\sum_{i=3}^n\mu_i^2+(n-\varepsilon)\sum_{k=1}^n\mu_k^2\\
&+\varepsilon\left(\tilde R_{1212}+\tilde R_{1234}\right)\\
&-\dfrac{\varepsilon}{2}\left[\sum_{i=1}^2\sum_{j,k=3}^n\mathring{\sigma}_{ijk}^2+\sum_{i,j=1}^2\sum_{k=3}^n\mathring{\sigma}_{ijk}^2+\sum_{k=1}^n\left(\sum_{i=1}^2\mathring{\sigma}_{iik}\right)^2\right]-\dfrac{1-\varepsilon}{n-2}\sum_{i=3}^n\sum_{j,k=1}^n\mathring{\sigma}_{ijk}^2\\
&-\varepsilon\sum_{k=5}^n\left(\mathring{\sigma}_{13k}\mathring{\sigma}_{24k}-\mathring{\sigma}_{14k}\mathring{\sigma}_{23k}\right)\\
&-\varepsilon\left[\left(\mathring{\sigma}_{113}-\mathring{\sigma}_{223}\right)\mathring{\sigma}_{124}-\left(\mathring{\sigma}_{114}-\mathring{\sigma}_{224}\right)\mathring{\sigma}_{123}+\left(\mathring{\sigma}_{331}-\mathring{\sigma}_{441}\right)\mathring{\sigma}_{234}-\left(\mathring{\sigma}_{332}-\mathring{\sigma}_{442}\right)\mathring{\sigma}_{134}\right]\\
\leq&\dfrac{(n-2)\varepsilon-2}{2}\sum_{k=1}^n\sum_{i=1}^2\mathring{\sigma}_{iik}\mu_k+\dfrac{(n-4)\varepsilon}{2}\sum_{i=1}^2\mu_i^2+(1-\varepsilon)\sum_{i=3}^n\mu_i^2+(n-\varepsilon)\sum_{k=1}^n\mu_k^2\\
&+\varepsilon\left(\tilde R_{1212}+\tilde R_{1234}\right)-\dfrac{3(1-\varepsilon)n}{2(n-2)^2}\sum_{k=3}^n\left(\sum_{i=1}^2\mathring{\sigma}_{iik}\right)^2\\
&-\dfrac{\varepsilon}{2}\left[\dfrac12\sum_{j=3}^n\left(\mathring{\sigma}_{11j}+\mathring{\sigma}_{22j}\right)^2+\dfrac12\sum_{i=1}^2\left(\mathring{\sigma}_{i33}+\mathring{\sigma}_{i44}\right)^2+\sum_{i=1}^2\sum_{j,k=5}^n\mathring{\sigma}_{ijk}^2+\sum_{k=1}^n\left(\sum_{i=1}^2\mathring{\sigma}_{iik}\right)^2\right]\\
\leq&\dfrac{(n-2)\varepsilon-2}{2}\sum_{k=1}^n\sum_{i=1}^2\mathring{\sigma}_{iik}\mu_k+\dfrac{(n-4)\varepsilon}{2}\sum_{i=1}^2\mu_i^2+(1-\varepsilon)\sum_{i=3}^n\mu_i^2+(n-\varepsilon)\sum_{k=1}^n\mu_k^2\\
&+\varepsilon\left(\tilde R_{1212}+\tilde R_{1234}\right)-\dfrac{3(1-\varepsilon)n}{2(n-2)^2}\sum_{k=3}^n\left(\sum_{i=1}^2\mathring{\sigma}_{iik}\right)^2\\
&-\dfrac{\varepsilon}{2}\left[\dfrac32\sum_{j=3}^n\left(\mathring{\sigma}_{11j}+\mathring{\sigma}_{22j}\right)^2+\dfrac{n-1}{n-2}\sum_{j=1}^2\left(\mathring{\sigma}_{11j}+\mathring{\sigma}_{22j}\right)^2\right]\\
\leq&\left[n-\varepsilon+\dfrac{(n-4)\varepsilon}{2}+\dfrac{(n-2)((n-2)\varepsilon-2)^2}{8(n-1)\varepsilon}\right]\sum_{i=1}^2\mu_j^2\\
&+\left[n+1-2\varepsilon+\dfrac{(n-2)^2((n-2)\varepsilon-2)^2}{12\left((n^2-6n+4)\varepsilon+2n\right)}\right]\sum_{j=3}^n\mu_j^2+\varepsilon\left(\tilde R_{1212}+\tilde R_{1234}\right).
\end{align*}

Taking $\varepsilon=\frac{2}{n-2}$, we get
\begin{align*}
\dfrac{1}{n-2}\sum_{i=1}^2\tilde R_{ii}+\dfrac{n-4}{(n-2)^2}\sum_{i=3}^n\tilde R_{ii}\leq&\dfrac{(n+2)(n-3)}{n-2}\sum_{k=1}^n\mu_k^2+\dfrac{2}{n-2}\left(\tilde R_{1212}+\tilde R_{1234}\right).
\end{align*}
\hfill Q.E.D.

\begin{lemma}\label{lemma3.7}Let $\sigma$ and $\tilde R$ be as in Lemma \ref{lemma3.4}. Then for all orthonormal frame $\{e_1,e_2,e_3,e_4,\cdots, e_n\}$,
\begin{align}\label{E-lemma3.7}
\sum_{i=1}^2\sum_{j=3}^4\tilde R_{ijij}-2\tilde R_{1234}\geq \frac{2}{3}\left[\sum_{i=1}^4\sum_{j=1}^n\tilde R_{ijij}-\frac{n-3}{3n-8}\sum_{k=1}^n\left(H^{k^*}\right)^2\right].
\end{align}
\end{lemma}

\noindent \textbf{Proof:} 
Using notations in Lemma \ref{lemma3.6}, we get
\begin{align}\label{EEE3.1}
\dfrac{1}{4}\sum_{i=1}^4\tilde R_{ii}=&n\sum_{k=1}^n\mu_k^2+\dfrac{n-2}{4}\left(\sum_{k=1}^n\sum_{i=1}^4\mathring{\sigma}_{iik}\mu_k+\sum_{i=1}^4\mu_i^2\right)-\dfrac{1}{4}\sum_{i=1}^4\sum_{j,k=1}^n\mathring{\sigma}_{ijk}^2.
\end{align}
By (\ref{e-1234}) and (\ref{e-cur}), we conclude that,
\begin{align}\label{EEE3.2}
\sum_{i=1}^2\sum_{j=3}^4\tilde{R}_{ijij}-\tilde R_{1234}\geq&\dfrac12\sum_{k=1}^n\left(\sum_{i=1}^4\mathring{\sigma}_{iik}\right)^2+2\sum_{k=1}^n\sum_{i=1}^4\mathring{\sigma}_{iik}\mu_k+2\sum_{i=1}^4\mu_i^2+4\sum_{k=1}^n\mu_k^2\nonumber\\
&-\sum_{k=1}^2\left(\mathring{\sigma}_{33k}+\mathring{\sigma}_{44k}\right)^2-\sum_{k=3}^4\left(\mathring{\sigma}_{11k}+\mathring{\sigma}_{22k}\right)^2-\dfrac12\sum_{k=5}^n\left(\mathring{\sigma}_{11k}+\mathring{\sigma}_{22k}\right)^2\nonumber\\
&-\dfrac12\sum_{k=1}^2\left(\mathring{\sigma}_{11k}+\mathring{\sigma}_{22k}\right)^2-\dfrac12\sum_{k=3}^4\left(\mathring{\sigma}_{33k}+\mathring{\sigma}_{44k}\right)^2-\dfrac12\sum_{k=5}^n\left(\mathring{\sigma}_{33k}+\mathring{\sigma}_{44k}\right)^2\nonumber\\
&-\sum_{k=3}^4\left(\mathring{\sigma}_{11k}-\mathring{\sigma}_{22k}\right)^2-\sum_{k=1}^2\left(\mathring{\sigma}_{33k}-\mathring{\sigma}_{44k}\right)^2\nonumber\\
&-4\sum_{k=3}^4\mathring{\sigma}_{12k}^2-4\sum_{i=1}^2\mathring{\sigma}_{i34}^2-2\sum_{i=1}^2\sum_{j=3}^4\sum_{k=5}^n\mathring{\sigma}_{ijk}^2.
\end{align}
On the other hand,
\begin{align}\label{EEE3.3}
\sum_{i=1}^4\sum_{j,k=1}^n\mathring{\sigma}_{ijk}^2
=&\sum_{i,j,k=1}^2\mathring{\sigma}_{ijk}^2+\sum_{i,j,k=3}^4\mathring{\sigma}_{ijk}^2+3\sum_{i=1}^2\sum_{j=3}^4\mathring{\sigma}_{ijj}^2+3\sum_{i=1}^2\sum_{k=3}^4\mathring{\sigma}_{iik}^2\nonumber\\
&+2\sum_{i=1}^4\sum_{k=5}^n\mathring{\sigma}_{iik}^2+\sum_{i=1}^4\sum_{j=5}^n\mathring{\sigma}_{ijj}^2\nonumber\\
&+6\sum_{i=1}^2\mathring{\sigma}_{i34}^2+6\sum_{k=3}^4\mathring{\sigma}_{12k}^2+4\sum_{1\leq i<j\leq 4}\sum_{k=5}^n\mathring{\sigma}_{ijk}^2+2\sum_{i=1}^4\sum_{5\leq j<k\leq n}\mathring{\sigma}_{ijk}^2\nonumber\\
\geq&\dfrac34\sum_{i=1}^2\left(\sum_{j=1}^2\mathring{\sigma}_{ijj}\right)^2+\dfrac34\sum_{i=3}^4\left(\sum_{j=3}^4\mathring{\sigma}_{ijj}\right)^2\nonumber\\
&+\dfrac32\sum_{i=1}^2\left(\mathring{\sigma}_{i33}+\mathring{\sigma}_{i44}\right)^2+\dfrac32\sum_{k=3}^4\left(\mathring{\sigma}_{11k}+\mathring{\sigma}_{22k}\right)^2+\dfrac32\sum_{i=1}^2\left(\mathring{\sigma}_{i33}-\mathring{\sigma}_{i44}\right)^2+\dfrac32\sum_{k=3}^4\left(\mathring{\sigma}_{11k}-\mathring{\sigma}_{22k}\right)^2\nonumber\\
&+\dfrac34\sum_{k=5}^n\left(\mathring{\sigma}_{11k}+\mathring{\sigma}_{22k}\right)^2+\dfrac34\sum_{k=5}^n\left(\mathring{\sigma}_{33k}+\mathring{\sigma}_{44k}\right)^2\nonumber\\
&+\dfrac{1}{8}\sum_{k=5}^n\left(\sum_{i=1}^4\mathring{\sigma}_{iik}\right)^2+\sum_{i=1}^4\sum_{j=5}^n\mathring{\sigma}_{ijj}^2+6\sum_{i=1}^2\mathring{\sigma}_{i34}^2+6\sum_{k=3}^4\mathring{\sigma}_{12k}^2+4\sum_{i=1}^2\sum_{j=3}^4\sum_{k=5}^n\mathring{\sigma}_{ijk}^2.
\end{align}

\vspace{.1in}

If $n=4$, we have 
\begin{equation*}
\sum_{i=1}^4\mathring{\sigma}_{iik}=0, \ \ \ for \ 1\leq k\leq 4.
\end{equation*}
From (\ref{EEE3.1}), (\ref{EEE3.2}) and (\ref{EEE3.3}), we compute
\begin{align}\label{EEE3.4}
\sum_{i=1}^2\sum_{j=3}^4\tilde{R}_{ijij}-\tilde R_{1234}
\geq&6\sum_{i=1}^4\mu_i^2-\sum_{k=1}^2\left(\mathring{\sigma}_{33k}+\mathring{\sigma}_{44k}\right)^2-\sum_{k=3}^4\left(\mathring{\sigma}_{11k}+\mathring{\sigma}_{22k}\right)^2-\sum_{k=3}^4\left(\mathring{\sigma}_{11k}-\mathring{\sigma}_{22k}\right)^2\nonumber\\
&-\sum_{k=1}^2\left(\mathring{\sigma}_{33k}-\mathring{\sigma}_{44k}\right)^2
-\dfrac12\sum_{k=1}^2\left(\mathring{\sigma}_{11k}+\mathring{\sigma}_{22k}\right)^2-\dfrac12\sum_{k=3}^4\left(\mathring{\sigma}_{33k}+\mathring{\sigma}_{44k}\right)^2\nonumber\\
&-4\sum_{k=3}^4\mathring{\sigma}_{12k}^2-4\sum_{i=1}^2\mathring{\sigma}_{i34}^2\nonumber\\
\geq& 6\sum_{i=1}^4\mu_i^2-\frac{2}{3}\sum_{i,j,k=1}^4\mathring{\sigma}_{ijk}^2\nonumber\\
=& \dfrac{2}{3}\sum_{i=1}^4\sum_{j=1}^n\tilde R_{ijij}-\dfrac{1}{6}\sum_{k=1}^n\left(H^{k^*}\right)^2.
\end{align}

\vspace{.1in}

If $n\geq 5$, we have from (\ref{EEE3.1}), (\ref{EEE3.2}) and (\ref{EEE3.3}) that
\begin{align}\label{EEE3.5}
\sum_{i=1}^2\sum_{j=3}^4\tilde{R}_{ijij}-\tilde R_{1234}
\geq&\dfrac12\sum_{k=1}^n\left(\sum_{i=1}^4\mathring{\sigma}_{iik}\right)^2+2\sum_{k=1}^n\sum_{i=1}^4\mathring{\sigma}_{iik}\mu_k+2\sum_{i=1}^4\mu_i^2+4\sum_{k=1}^n\mu_k^2\nonumber\\
 &-\frac{2}{3}\left[\sum_{i=1}^4\sum_{j,k=1}^n\mathring{\sigma}_{ijk}^2-\dfrac{1}{8}\sum_{k=5}^n\left(\sum_{i=1}^4\mathring{\sigma}_{iik}\right)^2-\sum_{i=1}^4\sum_{j=5}^n\mathring{\sigma}_{ijj}^2\right]\nonumber\\
\geq&\dfrac12\sum_{k=1}^n\left(\sum_{i=1}^4\mathring{\sigma}_{iik}\right)^2-\frac{2(n-5)}{3}\sum_{k=1}^n\sum_{i=1}^4\mathring{\sigma}_{iik}\mu_k\nonumber\\
& -\frac{2(n-5)}{3}\sum_{i=1}^4\mu_i^2-\left(\frac{8n}{3}-4\right)\sum_{k=1}^n\mu_k^2\nonumber\\
 &+\frac{2}{3}\sum_{i=1}^4\tilde R_{ii}+\dfrac{1}{12}\sum_{k=5}^n\left(\sum_{i=1}^4\mathring{\sigma}_{iik}\right)^2+\frac{2}{3(n-4)}\sum_{i=1}^4\left(\sum_{j=1}^4\mathring{\sigma}_{ijj}\right)^2.
  \end{align}

\vspace{.1in}

Then from (\ref{EEE3.5}), we have
\begin{align}\label{EEE3.7}
\sum_{i=1}^2\sum_{j=3}^4\tilde{R}_{ijij}-\tilde R_{1234}
\geq&\frac{2}{3}\sum_{i=1}^4\sum_{j=1}^n\tilde R_{ijij}+\frac{3n-8}{6(n-4)}\sum_{k=1}^4\left(\sum_{i=1}^4\mathring{\sigma}_{iik}\right)^2-\frac{2(n-5)}{3}\sum_{k=1}^4\sum_{i=1}^4\mathring{\sigma}_{iik}\mu_k\nonumber\\
& -\left(\frac{2(n-5)}{3}+\frac{8n}{3}-4\right)\sum_{k=1}^4\mu_k^2\nonumber\\
 &+\frac{7}{12}\sum_{k=5}^n\left(\sum_{i=1}^4\mathring{\sigma}_{iik}\right)^2-\frac{2(n-5)}{3}\sum_{k=5}^n\sum_{i=1}^4\mathring{\sigma}_{iik}\mu_k-\left(\frac{8n}{3}-4\right)\sum_{k=5}^n\mu_k^2\nonumber\\
\geq&\frac{2}{3}\left[\sum_{i=1}^4\sum_{j=1}^n\tilde R_{ijij}-\frac{(n-3)(n+2)^2}{3n-8}\sum_{k=1}^4\mu_k^2-\frac{2(n+2)^2}{7}\sum_{k=5}^n\mu_k^2\right]\nonumber\\
\geq&\frac{2}{3}\left[\sum_{i=1}^4\sum_{j=1}^n\tilde R_{ijij}-\frac{(n-3)(n+2)^2}{3n-8}\sum_{k=1}^n\mu_k^2\right]\nonumber\\
=&\frac{2}{3}\left[\sum_{i=1}^4\sum_{j=1}^n\tilde R_{ijij}-\frac{n-3}{3n-8}\sum_{k=1}^n\left(H^{k^*}\right)^2\right],
\end{align}
where we have used the fact that
\begin{equation*}
\frac{n-3}{3n-8}-\frac{2}{7}=\frac{n-5}{7(3n-8)}\geq 0
\end{equation*}
for $n\geq5$.

\vspace{.1in}

From (\ref{EEE3.4}) and (\ref{EEE3.7}), we can easily see that we have the unified estimation (\ref{E-lemma3.7}). This finishes the proof of the lemma.
\hfill Q.E.D.

\vspace{.2in}

\section{Sphere theorems for Lagrangian submanifolds in  K\"ahler manifold}

\vspace{.1in}

In this section, we will prove the sphere theorems for Lagrangian submanifolds in K\"ahler manifold. Let's first examine more about the curvature tensor on a Lagrangian submanifold.

\begin{proposition} Let $M^n$ be a Lagrangian submanifold of a K\"ahler manifold $N^{2n}$. Then for any orthonormal unit vector fields $X,Y,Z,W$ tangent to $M$, we have
 \begin{equation}\label{E-Lag1}
 \frac{3}{4}\tilde{K}_{\min}-\frac{1}{2}\tilde{K}_{\max}\leq K(X,Y)\leq \frac{3}{4}\tilde{K}_{\max}-\frac{1}{2}\tilde{K}_{\min},
\end{equation}
and
 \begin{equation}\label{E-Lag2}
 \frac{1}{2}(\tilde{K}_{\min}-\tilde{K}_{\max})\leq K(X,Y,Z,W)\leq \frac{1}{2}(\tilde{K}_{\max}-\tilde{K}_{\min}).
\end{equation}
\end{proposition}

\noindent \textbf{Proof:}
 By (\ref{e1.4}), we have for any vector field $X$ on $N$ that
 \begin{equation}\label{E4.3}
\tilde{K}_{\min}|X|^4\leq K(X)\leq \tilde{K}_{\max}|X|^4.
\end{equation}
By (\ref{E2.2}) and (\ref{E4.3}), we have for any orthonormal vector fields $X,Y$ on $N$
\begin{eqnarray*}
32K(X,Y)
&\leq& 3\tilde{K}_{\max}\left(|X+JY|^4+|X-JY|^4\right)\nonumber\\
&  & -\tilde{K}_{\min}\left(|X+Y|^4+|X-Y|^4+4|X|^4+4|Y|^4)\right)\nonumber\\
&=& 24(1+\langle X, JY\rangle^2)\tilde{K}_{\max}-16\tilde{K}_{\min}.
\end{eqnarray*}
Similarly we have
 \begin{equation*}
32K(X,Y)\geq 24(1+\langle X, JY\rangle^2)\tilde{K}_{\min}-16\tilde{K}_{\max}.
\end{equation*}
Therefore, we have
 \begin{equation*}
 \frac{3}{4}(1+\langle X, JY\rangle^2)\tilde{K}_{\min}-\frac{1}{2}\tilde{K}_{\max}\leq K(X,Y)\leq \frac{3}{4}(1+\langle X, JY\rangle^2)\tilde{K}_{\max}-\frac{1}{2}\tilde{K}_{\min}.
\end{equation*}
Since $M$ is Lagrangian, (\ref{E-Lag1}) follows.

By (\ref{E2.3}) and (\ref{E4.3}), we have for any orthonormal vector fields $X,Y,Z,W$ on $N$
\begin{eqnarray*}
256K(X,Y,Z,W)
&\leq& \tilde{K}_{\max}\left(|X+Z+JY+JW|^4+|X+Z-JY-JW|^4\right.\nonumber\\
&  & \ \ \ \ \ \ \left.  +|X-Z+JY-JW|^4+|X-Z-JY+JW|^4   \right.\nonumber\\
&  & \ \ \ \ \ \ \left.  +|X+W+JY-JZ|^4+|X+W-JY+JZ|^4   \right.\nonumber\\
&  & \ \ \ \ \ \ \left.  +|X-W+JY+JZ|^4+|X-W-JY-JZ|^4   \right)\nonumber\\
&  & -\tilde{K}_{\min}\left(|X+Z+JY-JW|^4+|X+Z-JY+JW|^4\right.\nonumber\\
&  & \ \ \ \ \ \ \ \ \left.  +|X-Z+JY+JW|^4+|X-Z-JY-JW|^4   \right.\nonumber\\
&  & \ \ \ \ \ \ \ \ \left.  +|X+W+JY+JZ|^4+|X+W-JY-JZ|^4   \right.\nonumber\\
&  & \ \ \ \ \ \ \ \ \left.  +|X-W+JY-JZ|^4+|X-W-JY+JZ|^4   \right)\nonumber\\
&=& \tilde{K}_{\max}\left[128+8(\langle X+Z,JY+JW\rangle^2+\langle X-Z,JY-JW\rangle^2  \right.\nonumber\\
&  & \ \ \ \ \ \  \ \ \ \ \ \ \ \ \ \ \ \left.  +\langle X+W,JY-JZ\rangle^2+\langle X-W,JY+JZ\rangle^2)   \right]\nonumber\\
&  & -\tilde{K}_{\min}\left[128+8(\langle X+Z,JY-JW\rangle^2+\langle X-Z,JY+JW\rangle^2  \right.\nonumber\\
&  & \ \ \ \ \ \  \ \ \ \ \ \ \ \ \ \ \ \ \ \left.  +\langle X+W,JY+JZ\rangle^2+\langle X-W,JY-JZ\rangle^2)   \right].
\end{eqnarray*}
Similarly, we have
\begin{eqnarray*}
256K(X,Y,Z,W)
&\geq& \tilde{K}_{\min}\left[128+8(\langle X+Z,JY+JW\rangle^2+\langle X-Z,JY-JW\rangle^2  \right.\nonumber\\
&  & \ \ \ \ \ \  \ \ \ \ \ \ \ \ \ \ \ \left.  +\langle X+W,JY-JZ\rangle^2+\langle X-W,JY+JZ\rangle^2)   \right]\nonumber\\
&  & -\tilde{K}_{\max}\left[128+8(\langle X+Z,JY-JW\rangle^2+\langle X-Z,JY+JW\rangle^2  \right.\nonumber\\
&  & \ \ \ \ \ \  \ \ \ \ \ \ \ \ \ \ \ \ \ \left.  +\langle X+W,JY+JZ\rangle^2+\langle X-W,JY-JZ\rangle^2)   \right].
\end{eqnarray*}
Since $M$ is Lagrangian, (\ref{E-Lag2}) follows.
\hfill Q.E.D.

\vspace{.1in}

\noindent\textbf{Proof of Theorem A:} We will show that under our assumption, $M\times{\mathbb R}^2$ has nonnegative isotropic curvature, i.e., (\ref{e-BS}) holds for all orthonormal four-frames $\{e_1,e_2,e_3,e_4\}$ and all $\lambda, \mu\in[0,1]$. For that purpose, we first extend the four-frame $\{e_1,e_2,e_3,e_4\}$ to be an adapted frame $\{e_1,\cdots,e_{2n}\}$ of $N$ such that $\{e_1,\cdots,e_{n}\}$ are tangent to $M$ and $\{e_{n+1}=Je_1,\cdots,e_{2n}=Je_n\}$ are normal to $M$.  The Gauss equation (\ref{e-gauss}) implies that
\begin{equation}\label{E4.4}
\tilde R(X,Y,Z,W):=R(X,Y,Z,W)-K(X,Y,Z,W)
\end{equation}
is an algebraic curvature. 

First note that (cf. 
Lemma 3.1 of  \cite{CS}, with $e_4$ replaced by $-e_4$)
\begin{eqnarray}\label{e-cs}
12 K_{1234}
&= & -4(K_{1212}+K_{3434})-2(K_{1313}+K_{1414}+K_{2323}+K_{2424})\\
& & +\left[K(e_1+e_3,e_2+e_4)+K(e_1-e_3,e_2-e_4)+K(e_2+e_3,e_1-e_4)+K(e_2-e_3,e_1+e_4)\right].\nonumber
\end{eqnarray}
Introduce $K_{ij}:=\sum_{k=1}^nK_{ikjk}$. By \eqref{E-Lag1}, we have for every $0<\varepsilon\leq 1$,
\begin{align*}
\dfrac{Ric^{[2]}_{\min}}{2}\leq&\dfrac{\varepsilon}{2}\dfrac{\left(R_{11}+R_{33}\right)+\lambda^2\left(R_{11}+R_{44}\right)+\mu^2\left(R_{22}+R_{33}\right)+\lambda^2\mu^2\left(R_{22}+R_{44}\right)}{\left(1+\lambda^2\right)\left(1+\mu^2\right)}\\
&+\dfrac{1-\varepsilon}{n-2}\left[\dfrac{\left(R_{22}+R_{44}\right)+\lambda^2\left(R_{22}+R_{33}\right)+\mu^2\left(R_{11}+R_{44}\right)+\lambda^2\mu^2\left(R_{11}+R_{33}\right)}{\left(1+\lambda^2\right)\left(1+\mu^2\right)}+\sum_{i=5}^nR_{ii}\right]\\
=&\dfrac{\varepsilon}{2}\dfrac{\left(K_{11}+K_{33}\right)+\lambda^2\left(K_{11}+K_{44}\right)+\mu^2\left(K_{22}+K_{33}\right)+\lambda^2\mu^2\left(K_{22}+K_{44}\right)}{\left(1+\lambda^2\right)\left(1+\mu^2\right)}\\
&+\dfrac{1-\varepsilon}{n-2}\left[\dfrac{\left(K_{22}+K_{44}\right)+\lambda^2\left(K_{22}+K_{33}\right)+\mu^2\left(K_{11}+K_{44}\right)+\lambda^2\mu^2\left(K_{11}+K_{33}\right)}{\left(1+\lambda^2\right)\left(1+\mu^2\right)}+\sum_{i=5}^nK_{ii}\right]\\
&+\dfrac{\varepsilon}{2}\dfrac{\left(\tilde R_{11}+\tilde R_{33}\right)+\lambda^2\left(\tilde R_{11}+\tilde R_{44}\right)+\mu^2\left(\tilde R_{22}+\tilde R_{33}\right)+\lambda^2\mu^2\left(\tilde R_{22}+\tilde R_{44}\right)}{\left(1+\lambda^2\right)\left(1+\mu^2\right)}\\
&+\dfrac{1-\varepsilon}{n-2}\left[\dfrac{\left(\tilde R_{22}+\tilde R_{44}\right)+\lambda^2\left(\tilde R_{22}+\tilde R_{33}\right)+\mu^2\left(\tilde R_{11}+\tilde R_{44}\right)+\lambda^2\mu^2\left(\tilde R_{11}+\tilde R_{33}\right)}{\left(1+\lambda^2\right)\left(1+\mu^2\right)}+\sum_{i=5}^n\tilde R_{ii}\right]\\
\leq&\dfrac{\varepsilon}{2}\dfrac{\left(K_{11}+K_{33}\right)+\lambda^2\left(K_{11}+K_{44}\right)+\mu^2\left(K_{22}+K_{33}\right)+\lambda^2\mu^2\left(K_{22}+K_{44}\right)}{\left(1+\lambda^2\right)\left(1+\mu^2\right)}\\
&-\varepsilon\dfrac{K_{1313}+\lambda^2K_{1414}+\mu^2K_{2323}+\lambda^2\mu^2K_{2424}-2\lambda\mu K_{1234}}{\left(1+\lambda^2\right)\left(1+\mu^2\right)}\\
&+\varepsilon\dfrac{R_{1313}+\lambda^2R_{1414}+\mu^2R_{2323}+\lambda^2\mu^2R_{2424}-2\lambda\mu R_{1234}}{\left(1+\lambda^2\right)\left(1+\mu^2\right)}\\
&-\varepsilon\dfrac{\tilde R_{1313}+\lambda^2\tilde R_{1414}+\mu^2\tilde R_{2323}+\lambda^2\mu^2\tilde R_{2424}-2\lambda\mu\tilde R_{1234}}{\left(1+\lambda^2\right)\left(1+\mu^2\right)}\\
&+(1-\varepsilon)(n-1)\left(\dfrac34\tilde K_{\max}-\dfrac12\tilde K_{\min}\right)\\
&+\dfrac{\varepsilon}{2}\dfrac{\left(\tilde R_{11}+\tilde R_{33}\right)+\lambda^2\left(\tilde R_{11}+\tilde R_{44}\right)+\mu^2\left(\tilde R_{22}+\tilde R_{33}\right)+\lambda^2\mu^2\left(\tilde R_{22}+\tilde R_{44}\right)}{\left(1+\lambda^2\right)\left(1+\mu^2\right)}\\
&+\dfrac{1-\varepsilon}{n-2}\left[\dfrac{\left(\tilde R_{22}+\tilde R_{44}\right)+\lambda^2\left(\tilde R_{22}+\tilde R_{33}\right)+\mu^2\left(\tilde R_{11}+\tilde R_{44}\right)+\lambda^2\mu^2\left(\tilde R_{11}+\tilde R_{33}\right)}{\left(1+\lambda^2\right)\left(1+\mu^2\right)}+\sum_{i=5}^n\tilde R_{ii}\right]\\
\leq&(n-1-\varepsilon)\left(\dfrac34\tilde K_{\max}-\dfrac12\tilde K_{\min}\right)\\
&+\varepsilon\dfrac{R_{1313}+\lambda^2R_{1414}+\mu^2R_{2323}+\lambda^2\mu^2R_{2424}-2\lambda\mu R_{1234}}{\left(1+\lambda^2\right)\left(1+\mu^2\right)}\\
&-\varepsilon\dfrac{\tilde R_{1313}+\lambda^2\tilde R_{1414}+\mu^2\tilde R_{2323}+\lambda^2\mu^2\tilde R_{2424}-2\lambda\mu\tilde R_{1234}}{\left(1+\lambda^2\right)\left(1+\mu^2\right)}\\
&+\dfrac{\varepsilon}{2}\dfrac{\left(\tilde R_{11}+\tilde R_{33}\right)+\lambda^2\left(\tilde R_{11}+\tilde R_{44}\right)+\mu^2\left(\tilde R_{22}+\tilde R_{33}\right)+\lambda^2\mu^2\left(\tilde R_{22}+\tilde R_{44}\right)}{\left(1+\lambda^2\right)\left(1+\mu^2\right)}\\
&+\dfrac{1-\varepsilon}{n-2}\left[\dfrac{\left(\tilde R_{22}+\tilde R_{44}\right)+\lambda^2\left(\tilde R_{22}+\tilde R_{33}\right)+\mu^2\left(\tilde R_{11}+\tilde R_{44}\right)+\lambda^2\mu^2\left(\tilde R_{11}+\tilde R_{33}\right)}{\left(1+\lambda^2\right)\left(1+\mu^2\right)}+\sum_{i=5}^n\tilde R_{ii}\right].
\end{align*}
Here we used the estimate from (\ref{E-Lag1}) and (\ref{e-cs}) that
\begin{align*}
&\dfrac{1}{2}\dfrac{\left(K_{11}+K_{33}\right)+\lambda^2\left(K_{11}+K_{44}\right)+\mu^2\left(K_{22}+K_{33}\right)+\lambda^2\mu^2\left(K_{22}+K_{44}\right)}{\left(1+\lambda^2\right)\left(1+\mu^2\right)}\\
&-\dfrac{K_{1313}+\lambda^2K_{1414}+\mu^2K_{2323}+\lambda^2\mu^2K_{2424}-2\lambda\mu K_{1234}}{\left(1+\lambda^2\right)\left(1+\mu^2\right)}\\
\leq &(n-4)\left(\dfrac34\tilde K_{\max}-\dfrac12\tilde K_{\min}\right) \\
& + \frac{1}{2}\sum_{j=1}^4\dfrac{\left(K_{1j1j}+K_{3j3j}\right)+\lambda^2\left(K_{1j1j}+K_{4j4j}\right)+\mu^2\left(K_{2j2j}+K_{3j3j}\right)+\lambda^2\mu^2\left(K_{2j2j}+K_{4j4j}\right)}{\left(1+\lambda^2\right)\left(1+\mu^2\right)}\\
&-\dfrac{K_{1313}+\lambda^2K_{1414}+\mu^2K_{2323}+\lambda^2\mu^2K_{2424}-2\lambda\mu K_{1234}}{\left(1+\lambda^2\right)\left(1+\mu^2\right)}\\
=&(n-4)\left(\dfrac34\tilde K_{\max}-\dfrac12\tilde K_{\min}\right)+\frac{1}{2}(K_{1212}+K_{3434}) \\
& +\frac{(\lambda^2+\mu^2)(K_{1313}+K_{2424})+(1+\lambda^2\mu^2)(K_{1414}+K_{2323})}{2\left(1+\lambda^2\right)\left(1+\mu^2\right)}+\frac{2\lambda\mu K_{1234}}{\left(1+\lambda^2\right)\left(1+\mu^2\right)} \\
=&(n-4)\left(\dfrac34\tilde K_{\max}-\dfrac12\tilde K_{\min}\right)+\left(\frac{1}{2}-\frac{2\lambda\mu}{3\left(1+\lambda^2\right)\left(1+\mu^2\right)} \right)(K_{1212}+K_{3434}) \\
& +\left(\frac{\lambda^2+\mu^2}{2\left(1+\lambda^2\right)\left(1+\mu^2\right)}-\frac{\lambda\mu}{3\left(1+\lambda^2\right)\left(1+\mu^2\right)} \right)(K_{1313}+K_{2424}) \\
& +\left(\frac{1+\lambda^2\mu^2}{2\left(1+\lambda^2\right)\left(1+\mu^2\right)}-\frac{\lambda\mu}{3\left(1+\lambda^2\right)\left(1+\mu^2\right)} \right)(K_{1414}+K_{2323}) \\
& +\frac{\lambda\mu}{6\left(1+\lambda^2\right)\left(1+\mu^2\right)}\left[K(e_1+e_3,e_2+e_4)+K(e_1-e_3,e_2-e_4)\right.\\
& \hspace{3.4cm} \left.+K(e_2+e_3,e_1-e_4)+K(e_2-e_3,e_1+e_4)\right] \\
\leq&(n-2)\left(\dfrac34\tilde K_{\max}-\dfrac12\tilde K_{\min}\right).
\end{align*}
Now applying Lemma \ref{lemma3.6}, by taking $\varepsilon=\frac{2}{n-2}$, we conclude that
\begin{align*}
\dfrac{Ric_{\min}^{[2]}}{2}\leq&\left(n-1-\dfrac{2}{n-2}\right)\left(\dfrac34\tilde K_{\max}-\dfrac12\tilde K_{\min}\right)\\
&+\dfrac{2}{n-2}\dfrac{R_{1313}+\lambda^2R_{1414}+\mu^2R_{2323}+\lambda^2\mu^2R_{2424}-2\lambda\mu R_{1234}}{\left(1+\lambda^2\right)\left(1+\mu^2\right)}\\
&+\dfrac{n-3}{(n+2)(n-2)}|\bf{H}|^2.
\end{align*}
Thus
\begin{eqnarray*}
& & \frac{4}{n-2}(R_{1313}+\lambda^2R_{1414}+\mu^2R_{2323}+\lambda^2\mu^2R_{2424}-2\lambda\mu R_{1234})\nonumber\\
&\geq& (1+\lambda^2)(1+\mu^2)\left[Ric^{[2]}_{\min}-\dfrac{n(n-3)}{n-2}\left(\dfrac32\tilde K_{\max}-\tilde K_{\min}\right)-\dfrac{2(n-3)}{(n-2)(n+2)}|{\bf H}|^2\right]\nonumber\\
&\geq&0,
\end{eqnarray*}
the strict inequality holds for some point $x_0\in M$, where the last inequality follows from our assumption (\ref{e-A}). Hence $M$ is diffeomorphic to ${\mathbb S}^n$ by standard argument using Lemma \ref{lemma-Aubin}, Lemma \ref{lemma-MW}, Lemma \ref{lemma-Se} and Theorem \ref{thmBS} (see, for example, \cite{GX}, \cite{SS}).
\hfill Q.E.D.

\vspace{.1in}

Next, we turn to prove the topological sphere theorems for Lagrangian submanifolds. 

\vspace{.1in}

\noindent \textbf{Proof of Theorem B:} Define the operator $\tilde R$ by (\ref{E4.4}), which is an algebraic curvature. Then we have from Lemma \ref{lemma3.5} that for all orthonormal frame $\{e_1,e_2,e_3,e_4,\cdots, e_n\}$,
\begin{align*}
\sum_{i=1}^2\sum_{j=3}^4\tilde R_{ijij}-2\tilde R_{1234}\geq&\tilde\eta(n)|{\bf H}|^2-\dfrac{2}{3}\sum_{i,j,k=1}^n\sigma_{ijk}^2=\tilde\eta(n)|{\bf H}|^2-\dfrac{2}{3}|{\bf B}|^2.
\end{align*}
In other word, by using (\ref{e-scalar}), we have
\begin{eqnarray}\label{E4.9}
& & R_{1313}+R_{1414}+R_{2323}+R_{2424}-2R_{1234}\nonumber\\
&\geq& K_{1313}+K_{1414}+K_{2323}+K_{2424}-2K_{1234} +\tilde\eta(n)|{\bf H}|^2-\dfrac{2}{3}|{\bf B}|^2\nonumber\\
& = & K_{1313}+K_{1414}+K_{2323}+K_{2424}-2K_{1234}-\frac{2}{3}\sum_{i,j=1}^nK_{ijij}+\dfrac{2}{3}R_M+\left(\tilde\eta(n)-\frac{2}{3}\right)|{\bf H}|^2.
\end{eqnarray}
It suffices to estimate the terms involving the curvature tensor on $N$. We will follow the argument as in the proof of Theorem 3.2 in \cite{CS}. 
For that purpose, we have from (\ref{e-cs})
\begin{eqnarray}\label{e4.8}
\sum_{1\leq i<j\leq 4}K_{ijij}& = &\dfrac{1}{8}\left[K(e_1+e_3,e_2+e_4)+K(e_1-e_3,e_2-e_4)+K(e_2+e_3,e_1-e_4)+K(e_2-e_3,e_1+e_4)\right]\nonumber\\
& & +\dfrac12\left(K_{1212}+K_{3434}\right)+\dfrac34\left(K_{1313}+K_{1414}+K_{2323}+K_{2424}-2K_{1234}\right).
\end{eqnarray}
On the other hand,
\begin{equation*}
\sum_{i,j=1}^nK_{ijij}= \sum_{i,j=5}^nK_{ijij}+2\sum_{i=1}^4\sum_{j=5}^nK_{ijij}+2\sum_{1\leq i<j\leq 4}K_{ijij}.
\end{equation*}
Hence,  using (\ref{E-Lag1}), we estimate
\begin{eqnarray}\label{E4.10}
& & \frac{3}{4}\left(K_{1313}+K_{1414}+K_{2323}+K_{2424}-2K_{1234}-\frac{2}{3}\sum_{i,j=1}^nK_{ijij}\right)\nonumber\\
& = & -\frac{1}{2}\sum_{i,j=5}^nK_{ijij}-\sum_{i=1}^4\sum_{j=5}^nK_{ijij}-\frac{1}{2}(K_{1212}+K_{3434})\nonumber\\
& & -\frac{1}{8}\left[K(e_1+e_3,e_2+e_4)+K(e_1-e_3,e_2-e_4)+K(e_2+e_3,e_1-e_4)+K(e_2-e_3,e_1+e_4)\right]\nonumber\\
&\geq& -\frac{1}{2}(n-4)(n-5)\left(\frac{3}{4}\tilde{K}_{\max}-\frac{1}{2}\tilde{K}_{\min}\right)-4(n-4)\left(\frac{3}{4}\tilde{K}_{\max}-\frac{1}{2}\tilde{K}_{\min}\right)-\left(\frac{3}{4}\tilde{K}_{\max}-\frac{1}{2}\tilde{K}_{\min}\right)\nonumber\\
&  & -\frac{1}{8}\cdot4\left(\frac{3}{4}\tilde{K}_{\max}-\frac{1}{2}\tilde{K}_{\min}\right)\cdot4\nonumber\\
&=& -\frac{n^2-n-6}{2}\left(\frac{3}{4}\tilde{K}_{\max}-\frac{1}{2}\tilde{K}_{\min}\right).
\end{eqnarray}
Inserting (\ref{E4.10}) into (\ref{E4.9}), we have
\begin{eqnarray*}
&  & R_{1313}+R_{1414}+R_{2323}+R_{2424}-2R_{1234}\nonumber\\
&\geq& \frac{2}{3}\left[R_M-\frac{(n-3)(n+2)}{4}(3\tilde{K}_{\max}-2\tilde{K}_{\min})-\eta(n)|{\bf H}|^2\right]\nonumber\\
&\geq & 0,
\end{eqnarray*}
the strict inequality holds for some point $x_0\in M$, where the last inequality follows from our assumption (\ref{e-B}). Here, $\eta(n)$ is given by
\begin{equation*}
\eta(n)=1-\frac{3}{2}\tilde\eta(n)=
       \begin{cases}
          \hspace{0.53cm}      \frac{1}{4}, \ \ if \ n=4,\\
               \frac{3n-7}{3n+2}, \ \ if \ n\geq 5.
       \end{cases}
\end{equation*}
By Lemma \ref{lemma-Se}, $M$ admits a metric with positive isotropic curvature. Since $M$ is simply connected, $M$ is homeomorphic to ${\mathbb S}^n$ by Lemma \ref{lemma-MM}.
\hfill Q.E.D.

\vspace{.1in}

\noindent \textbf{Proof of Theorem C:}  Using the same notations as in the proof of Theorem B, we have
\begin{equation*}
\sum_{i=1}^4\sum_{j=1}^n\tilde R_{ijij}
= \sum_{i=1}^4Ric_{ii}-\sum_{i=1}^4\sum_{j=1}^nK_{ijij}.
\end{equation*}
By Lemma \ref{lemma3.7}, we obtain for all orthonormal frame $\{e_1,e_2,e_3,e_4,\cdots, e_n\}$,
\begin{equation*}
\sum_{i=1}^2\sum_{j=3}^4\tilde R_{ijij}-2\tilde R_{1234}\geq \frac{2}{3}\left[Ric^{[4]}_{\min}-\sum_{i=1}^4\sum_{j=1}^nK_{ijij}-\frac{n-3}{3n-8}|{\bf H}|^2\right].
\end{equation*}
In other word,
\begin{eqnarray}\label{E4.12}
& & \frac{3}{2}(R_{1313}+R_{1414}+R_{2323}+R_{2424}-2R_{1234})\nonumber\\
&\geq& \frac{3}{2}(K_{1313}+K_{1414}+K_{2323}+K_{2424}-2K_{1234})+Ric^{[4]}_{\min}-\sum_{i=1}^4\sum_{j=1}^nK_{ijij}-\frac{n-3}{3n-8}|{\bf H}|^2.
\end{eqnarray}
We need to estimate the terms involving the curvature tensor on $N$. As in the proof of Theorem B, by (\ref{e4.8}), we  have
\begin{eqnarray*}
\sum_{i=1}^4\sum_{j=1}^nK_{ijij}&= & 2\sum_{1\leq i<j\leq 4}K_{ijij}+\sum_{i=1}^4\sum_{j=5}^nK_{ijij}\nonumber\\
& = &\dfrac{1}{4}\left[K(e_1+e_3,e_2+e_4)+K(e_1-e_3,e_2-e_4)+K(e_2+e_3,e_1-e_4)+K(e_2-e_3,e_1+e_4)\right]\nonumber\\
& &+ \left(K_{1212}+K_{3434}\right) +\dfrac32\left(K_{1313}+K_{1414}+K_{2323}+K_{2424}-2K_{1234}\right)+\sum_{i=1}^4\sum_{j=5}^nK_{ijij}.
\end{eqnarray*}
Use (\ref{E-Lag1}) to estimate
\begin{eqnarray}\label{E4.13}
& & \frac{3}{2}\left(K_{1313}+K_{1414}+K_{2323}+K_{2424}-2K_{1234}\right)-\sum_{i=1}^4\sum_{j=1}^nK_{ijij}\nonumber\\
& = & -\sum_{i=1}^4\sum_{j=5}^nK_{ijij}-(K_{1212}+K_{3434})\nonumber\\
& & -\frac{1}{4}\left[K(e_1+e_3,e_2+e_4)+K(e_1-e_3,e_2-e_4)+K(e_2+e_3,e_1-e_4)+K(e_2-e_3,e_1+e_4)\right]\nonumber\\
&\geq& -4(n-4)\left(\frac{3}{4}\tilde{K}_{\max}-\frac{1}{2}\tilde{K}_{\min}\right)-2\left(\frac{3}{4}\tilde{K}_{\max}-\frac{1}{2}\tilde{K}_{\min}\right) -\frac{1}{4}\cdot4\left(\frac{3}{4}\tilde{K}_{\max}-\frac{1}{2}\tilde{K}_{\min}\right)\cdot4\nonumber\\
&=& -(2n-5)\left(\frac{3}{2}\tilde{K}_{\max}-\tilde{K}_{\min}\right).
\end{eqnarray}
Inserting (\ref{E4.13}) into (\ref{E4.12}), we have
\begin{eqnarray*}
&  & R_{1313}+R_{1414}+R_{2323}+R_{2424}-2R_{1234}\nonumber\\
&\geq& \frac{2}{3}\left[Ric^{[4]}_{\min}-(2n-5)\left(\frac{3}{2}\tilde{K}_{\max}-\tilde{K}_{\min}\right)-\frac{n-3}{3n-8}|{\bf H}|^2\right]\nonumber\\
&\geq & 0,
\end{eqnarray*}
the strict inequality holds for some point $x_0\in M$, where the last inequality follows from our assumption (\ref{e-C}). The theorem then follows.
\hfill Q.E.D.

\vspace{.2in}

\section{Sphere theorems for Legendrian submanifolds in  Sasaki space form}

\vspace{.1in}

In this section, we will prove the sphere theorems for Legendrian submanifolds in Sasaki space form. The following lemma is an easy consequence of Proposition \ref{prop2.3} and the definition of Legendrian submanifold: 

\begin{proposition} Let $M^n$ be a Legendrian submanifold of a Sasaki space form $N^{2n+1}(c)$, then for any vector fields $X,Y,Z,W$ tangent to $M$, we have
 \begin{equation*}
K(X,Y)W= \frac{1}{4}(c+3)\left(\langle Y,W\rangle X-\langle X,W\rangle Y\right),
\end{equation*}
and
 \begin{equation}\label{E-Leg2}
K(X,Y,Z,W)= \frac{1}{4}(c+3)\left(\langle X,Z\rangle \langle Y,W\rangle-\langle X,W\rangle \langle Y,Z\rangle\right).
\end{equation}
\end{proposition}

\vspace{.1in}

The following differentiable sphere can be viewed as a Legendrian correspondence of Theorem 1.1':

\vspace{.1in}

\begin{theorem}
Let $M$ be a smooth $n(\geq 4)$-dimensional closed simply connected Legendrian submanifold of  a  Sasaki space form $N^{2n+1}(c)$. If $M$ satisfies the following condition:
\begin{equation}\label{e-D}
R_M \geq \frac{(n-2)(n+1)}{4}(c+3)+\frac{2n-3}{2n+3}|{\bf H}|^2,
\end{equation}
and the strict inequality holds for some point $x_0\in M$. Then $M$ is diffeomorphic to ${\mathbb S}^n$.
\end{theorem}

\vspace{.1in}

\noindent \textbf{Proof:} For any orthonormal four-frame $\{e_1,e_2,e_3,e_4\}$, we extend it to be an adapted orthonormal frame $\{e_1,\cdots,e_{2n+1}\}$ of $N$ such that $\{e_1,\cdots,e_{n}\}$ are tangent to $M$ and $\{e_{n+1}=\phi e_1,\cdots,e_{2n}=\phi e_n, e_{2n+1}=\xi\}$ are normal to $M$.  The Gauss equation (\ref{e-gauss}) implies that
\begin{equation*}
\tilde R(X,Y,Z,W):=R(X,Y,Z,W)-K(X,Y,Z,W)
\end{equation*}
is an algebraic curvature. By Lemma \ref{lemma3.4}, (\ref{e5.9}) and (\ref{e-scalar}), we have
\begin{eqnarray*}
\tilde R_{1212}+\tilde R_{1234}
&\geq& \dfrac12\left(\dfrac{6}{2n+3}|{\bf H}|^2-|{\bf B}|^2\right)
=\frac{1}{2}\left(R_M-\sum_{i,j=1}^nK_{ijij}-\dfrac{2n-3}{2n+3}|{\bf H}|^2\right).
\end{eqnarray*}
Lemma \ref{lemma3.1} implies that
\begin{eqnarray*}
&  & \tilde R_{1313}+\lambda^2\tilde R_{1414}+\mu^2\tilde R_{2323}+\lambda^2\mu^2\tilde R_{2424}-2\lambda\mu \tilde R_{1234}\nonumber\\
&\geq &\frac{(1+\lambda^2)(1+\mu^2)}{2}\left(R_M-\sum_{i,j=1}^nK_{ijij}-\dfrac{2n-3}{2n+3}|{\bf H}|^2\right),
\end{eqnarray*}
i.e.,
\begin{eqnarray}\label{E5.5}
& & 2(R_{1313}+\lambda^2R_{1414}+\mu^2R_{2323}+\lambda^2\mu^2R_{2424}-2\lambda\mu R_{1234})\nonumber\\
&\geq& 2(K_{1313}+\lambda^2K_{1414}+\mu^2K_{2323}+\lambda^2\mu^2K_{2424}-2\lambda\mu K_{1234})\nonumber\\
&  &+(1+\lambda^2)(1+\mu^2)\left(R_M-\sum_{i,j=1}^nK_{ijij}-\dfrac{2n-3}{2n+3}|{\bf H}|^2\right).
\end{eqnarray}
By (\ref{E-Leg2}),
\begin{equation*}
K_{ijij}=\frac{c+3}{4}, \ \forall i\neq j, \ \  \ \  K_{1234}=0.
\end{equation*}
Therefore,
\begin{eqnarray}\label{E5.6}
& & 2(K_{1313}+\lambda^2K_{1414}+\mu^2K_{2323}+\lambda^2\mu^2K_{2424}-2\lambda\mu K_{1234})-(1+\lambda^2)(1+\mu^2)\sum_{i,j=1}^nK_{ijij}\nonumber\\
&=& 2(1+\lambda^2)(1+\mu^2)\frac{c+3}{4}-(1+\lambda^2)(1+\mu^2)n(n-1)\frac{c+3}{4}\nonumber\\
&=& -(1+\lambda^2)(1+\mu^2)\frac{(n-2)(n+1)}{4}(c+3).
\end{eqnarray}
Inserting (\ref{E5.6}) into (\ref{E5.5}), we have
\begin{eqnarray*}
& & 2(R_{1313}+\lambda^2R_{1414}+\mu^2R_{2323}+\lambda^2\mu^2R_{2424}-2\lambda\mu R_{1234})\nonumber\\
&\geq& (1+\lambda^2)(1+\mu^2)\left[R_M-\frac{(n-2)(n+1)}{4}(c+3)-\dfrac{2n-3}{2n+3}|{\bf H}|^2\right]\nonumber\\
&\geq&0,
\end{eqnarray*}
the strict inequality holds for some point $x_0\in M$, where the last inequality follows from our assumption (\ref{e-D}). Then the theorem follows.
\hfill Q.E.D.

\vspace{.1in}

Similar to proof of Theorem A, Theorem B, Theorem C, with (\ref{E-Lag1}) and (\ref{E-Lag2}) replaced by (\ref{E-Leg2}), we can obtain the following sphere theorems for submanifolds in Sasaki space form under various curvature assumptions. Since the proofs are similar, we omit the details.

\vspace{.1in}

\begin{theorem}
Let $M$ be a smooth $n(\geq 4)$-dimensional closed simply connected Legendrian submanifold of  a  Sasaki space form $N^{2n+1}(c)$. If $M$ satisfies the following condition:
\begin{equation*}
\frac{Ric_{\min}^{[2]}}{2}\geq\frac{n(n-3)}{4(n-2)}(c+3)+\frac{n-3}{(n+2)(n-2)}|{\bf H}|^2,
\end{equation*}
and the strict inequality holds for some point $x_0\in M$. Then $M$ is diffeomorphic to ${\mathbb S}^n$.
\end{theorem}

\begin{theorem}
Let $M$ be a smooth $n(\geq 4)$-dimensional closed simply connected Legendrian submanifold of  a  Sasaki space form $N^{2n+1}(c)$. If $M$ satisfies the following condition:
\begin{equation*}
R_M \geq \frac{(n-3)(n+2)}{4}(c+3)+\eta(n)|{\bf H}|^2,
\end{equation*}
and the strict inequality holds for some point $x_0\in M$ where $\eta(n)$ is given by (\ref{E-eta}), then $M$ is homeomorphic to ${\mathbb S}^n$.
\end{theorem}

\begin{theorem}
Let $M$ be a smooth $n(\geq 4)$-dimensional closed simply connected Legendrian submanifold of  a  Sasaki space form $N^{2n+1}(c)$. If $M$ satisfies the following condition:
\begin{equation*}
Ric^{[4]}_{\min}\geq \frac{2n-5}{2}(c+3)+\frac{n-3}{3n-8}|{\bf H}|^2,
\end{equation*}
and the strict inequality holds for some point $x_0\in M$, then $M$ is homeomorphic to ${\mathbb S}^n$.
\end{theorem}

\vspace{.2in}


\end{document}